\input amstex
\let\myfrac=\frac
\input eplain
\let\frac=\myfrac
\input epsf




\loadeufm \loadmsam \loadmsbm
\message{symbol names}\UseAMSsymbols\message{,}

\font\myfontdefault=cmr10

\font\mytdmchapfont=cmb10 at 14pt
\font\mytdmheadfont=cmb10 at 10pt
\font\mytdmsubheadfont=cmr10
	
\magnification 1200
\newif\ifinappendices
\newif\ifundefinedreferences
\newif\ifchangedreferences
\newif\ifloadreferences
\newif\ifmakebiblio
\newif\ifmaketdm

\undefinedreferencesfalse
\changedreferencesfalse


\loadreferencestrue
\makebibliofalse
\maketdmfalse

\def\headpenalty{-400}     
\def\proclaimpenalty{-200} 

%
%

\def\alphanum#1{\ifcase #1 _\or A\or B\or C\or D\or E\or F\or G\or H\or I\or J\or K\or L\or M\or N\or O\or P\or Q\or R\or S\or T\or U\or V\or W\or X\or Y\or Z\fi}
\def\gobbleeight#1#2#3#4#5#6#7#8{}

\newwrite\references
\newwrite\tdm
\newwrite\biblio

\newcount\chapno
\newcount\headno
\newcount\subheadno
\newcount\procno
\newcount\figno
\newcount\citationno

\def\setcatcodes{%
\catcode`\!=0 \catcode`\\=11}%

\ifloadreferences
    {\catcode`\@=11 \catcode`\_=11%
    \global\def\_@citation@AndBarbBegZegh{1}
\global\def\_@citation@Aubin{2}
\global\def\_@citation@CaffNirSpr{3}
\global\def\_@citation@ChenGreen{4}
\global\def\_@citation@EpsMard{5}
\global\def\_@citation@Fried{6}
\global\def\_@citation@GuanSpruck{7}
\global\def\_@citation@HarveyLawson{8}
\global\def\_@citation@KamiA{9}
\global\def\_@citation@KamTan{10}
\global\def\_@citation@KapB{11}
\global\def\_@citation@KraSch{12}
\global\def\_@citation@Kulkarni{13}
\global\def\_@citation@LabA{14}
\global\def\_@citation@LabB{15}
\global\def\_@citation@MazzPac{16}
\global\def\_@citation@SmiA{17}
\global\def\_@citation@SmiNLD{18}
\global\def\_@citation@SmiH{19}
\global\def\_@citation@Thurston{20}
\global\def\_@proc@TheoremExistenceAndUniqueness{1.1}
\global\def\_@proc@TheoremFoliations{1.2}
\global\def\_@proc@TheoremFoliationsInQuasiFuchsian{1.3}
\global\def\_@proc@TheoremContinuousDependance{1.4}
\global\def\_@proc@TheoremFourDSpecialCase{1.6}
\global\def\_@proc@TheoremStructureOfKPEnd{1.7}
\global\def\_@head@SpecialLagrangianCurvature{3}
\global\def\_@head@HyperbolicEnds{4}
\global\def\_@head@ConformalStructures{5}
\global\def\_@proc@ThmKulkarniPinkallMetric{5.2}
\global\def\_@proc@PropMaximulBalls{5.3}
\global\def\_@head@HeadGBC{6}
\global\def\_@proc@LemmaStructureFoliation{6.1}
\global\def\_@proc@LemmaStructureQuotient{6.2}
\global\def\_@proc@LemmaStructureConstruction{6.3}
\global\def\_@proc@DefinitionGBC{6.4}
\global\def\_@proc@LemmaGBCGivesKP{6.5}
\global\def\_@proc@LemmaGBCMeansMaximality{6.6}
\global\def\_@proc@CorContainedInKPEnd{6.7}
\global\def\_@proc@LemmaKPIsUniqueMaximum{6.9}
\global\def\_@head@DerivativeOfSLOp{7}
\global\def\_@proc@DerivativeOfShapeOperator{7.1}
\global\def\_@proc@DerivativeOfSLIsLaplacian{7.2}
\global\def\_@proc@UsefulDecreasingFunction{7.3}
\global\def\_@proc@ZerothOrderCoeffIsPos{7.4}
\global\def\_@proc@LemmaInvertibilityOfDerivatives{7.5}
\global\def\_@head@DeformingEquivariantImmersions{8}
\global\def\_@proc@EquivariantDeformation{8.1}
\global\def\_@proc@CompactnessA{9.1}
\global\def\_@proc@CompactnessB{9.2}
\global\def\_@head@UpperAndLowerBounds{10}
\global\def\_@proc@DefinitionWeakCurvatureBounds{10.1}
\global\def\_@proc@LemmaGeomMaxPrinc{10.4}
\global\def\_@proc@LemmaLowerCurvatureBound{10.5}
\global\def\_@proc@CurvatureOfLevelSetsI{10.6}
\global\def\_@proc@LemmaUpperCurvatureBound{10.7}
\global\def\_@proc@CurvatureOfLevelSetsII{10.8}
\global\def\_@proc@UpperBoundInHyperbolicEnd{10.9}
\global\def\_@proc@LemmaUnifEquiv{11.1}
\global\def\_@proc@LemmaUniqueness{11.2}
\global\def\_@proc@TheoremLocalExistence{12.1}
\global\def\_@proc@LemmaExistence{12.2}
\global\def\_@head@QuasiFuchsianManifolds{13}
\global\def\_@proc@LemmaUpperBoundInFuchsianMfd{13.1}
\global\def\_@proc@ChenGreenberg{A.1}
\global\def\_@proc@ThmKamiCoveringMap{A.2}
\global\def\_@proc@LocationOfBoundary{A.3}
    }%
\else
    \openout\references=references.tex
\fi

\newcount\newchapflag 
\newcount\showpagenumflag 

\global\chapno = -1 
\global\citationno=0
\global\headno = 0
\global\subheadno = 0
\global\procno = 0
\global\figno = 0

\def\resetcounters{%
\global\headno = 0%
\global\subheadno = 0%
\global\procno = 0%
\global\figno = 0%
}

\global\newchapflag=0 
\global\showpagenumflag=0 

\def\chinfo{\ifinappendices\alphanum\chapno\else\the\chapno\fi}%
\def\headinfo{\ifinappendices\alphanum\headno\else\the\headno\fi}%
\def\subheadinfo{\headinfo.\the\subheadno}
\def\procinfo{\headinfo.\the\procno}
\def\figinfo{\headinfo.\the\figno}
\def\citationinfo{\the\citationno}%
\def\nextheadno{\global\advance\headno by 1 \global\subheadno = 0 \global\procno = 0}
\def\nextsubheadno{\global\advance\subheadno by 1}
\def\nextprocno{\global\advance\procno by 1 \procinfo}
\def\nextfigno{\global\advance\figno by 1 \figinfo}

{\global\let\noe=\noexpand%
%
%
\catcode`\@=11%
\catcode`\_=11%
\setcatcodes%
!global!def!_@@internal@@makeref#1{%
!global!expandafter!def!csname #1ref!endcsname##1{%
!csname _@#1@##1!endcsname%
!expandafter!ifx!csname _@#1@##1!endcsname!relax%
    !write16{#1 ##1 not defined - run saving references}%
    !undefinedreferencestrue%
!fi}}%
!global!def!_@@internal@@makelabel#1{%
!global!expandafter!def!csname #1label!endcsname##1{%
!edef!temptoken{!csname #1info!endcsname}%
!ifloadreferences%
    !expandafter!ifx!csname _@#1@##1!endcsname!relax%
        !write16{#1 ##1 not hitherto defined - rerun saving references}%
        !changedreferencestrue%
    !else%
        !expandafter!ifx!csname _@#1@##1!endcsname!temptoken%
        !else
            !write16{#1 ##1 reference has changed - rerun saving references}%
            !changedreferencestrue%
        !fi%
    !fi%
!else%
    !expandafter!edef!csname _@#1@##1!endcsname{!temptoken}%
    !edef!textoutput{!write!references{\global\def\_@#1@##1{!temptoken}}}%
    !textoutput%
!fi}}%
!global!def!makecounter#1{!_@@internal@@makelabel{#1}!_@@internal@@makeref{#1}}%
!unsetcatcodes%
}
\makecounter{ch}%
\makecounter{head}%
\makecounter{subhead}%
\makecounter{proc}%
\makecounter{fig}%
\makecounter{citation}%
\def\newref#1#2{%
\def\temptext{#2}%
\edef\bibliotextoutput{\expandafter\gobbleeight\meaning\temptext}%
\global\advance\citationno by 1\citationlabel{#1}%
\ifmakebiblio%
    \edef\fileoutput{\write\biblio{\noindent\hbox to 0pt{\hss$[\the\citationno]$}\hskip 0.2em\bibliotextoutput\medskip}}%
    \fileoutput%
\fi}%
\def\cite#1{%
$[\citationref{#1}]$%
\ifmakebiblio%
    \edef\fileoutput{\write\biblio{#1}}%
    \fileoutput%
\fi%
}%
%
%
%

\let\mypar=\par


\def\raggedleft{\leftskip=0pt plus 1fil \parfillskip=0pt}


\font\lettrinefont=cmr10 at 28pt
\def\lettrine #1[#2][#3]#4%
{\hangafter -#1 \hangindent #2
\noindent\hskip -#2 \vtop to 0pt{
\kern #3 \hbox to #2 {\lettrinefont #4\hss}\vss}}

\font\mylettrinefont=cmr10 at 28pt
\def\mylettrine #1[#2][#3][#4]#5%
{\hangafter -#1 \hangindent #2
\noindent\hskip -#2 \vtop to 0pt{
\kern #3 \hbox to #2 {\mylettrinefont #5\hss}\vss}}


\edef\Pagetitle={Blank}

\headline={\hfil\Pagetitle\hfil}

\footline={{\hfil\myfontdefault\folio\hfil}}

\def\nextoddpage
{
\newpage%
\ifodd\pageno%
\else%
    \global\showpagenumflag = 0%
    \null%
    \vfil%
    \eject%
    \global\showpagenumflag = 1%
\fi%
}


\def\newchap#1#2%
{%
%
%
\global\advance\chapno by 1%
\resetcounters%
%
%
\newpage%
\ifodd\pageno%
\else%
    \global\showpagenumflag = 0%
    \null%
    \vfil%
    \eject%
    \global\showpagenumflag = 1%
\fi%
\global\newchapflag = 1%
\global\showpagenumflag = 1%
%
%
{\font\chapfontA=cmsl10 at 30pt%
\font\chapfontB=cmsl10 at 25pt%
\null\vskip 5cm%
{\chapfontA\raggedleft\hfil%
{%
\ifnum\chapno=0
    \phantom{%
    \ifinappendices%
        Annexe \alphanum\chapno%
    \else%
        \the\chapno%
    \fi}%
\else%
    \ifinappendices%
        Annexe \alphanum\chapno%
    \else%
        \the\chapno%
    \fi%
\fi%
}%
\par}%
\vskip 2cm%
{\chapfontB\raggedleft%
\lineskiplimit=0pt%
\lineskip=0.8ex%
\hfil #1\par}%
\vskip 2cm%
}%
\edef\Pagetitle{#2}%
%
%
\ifmaketdm%
    \def\temp{#2}%
    \def\tempbis{\nobreak}%
    \edef\chaptitle{\expandafter\gobbleeight\meaning\temp}%
    \edef\mynobreak{\expandafter\gobbleeight\meaning\tempbis}%
    \edef\textoutput{\write\tdm{\bigskip{\noexpand\mytdmchapfont\noindent\chinfo\ - \chaptitle\hfill\noexpand\folio}\par\mynobreak}}%
\fi%
\textoutput%
}


\def\newhead#1%
{%
\ifhmode%
    \mypar%
\fi%
\ifnum\headno=0%
\ifinappendices
    \nobreak\vskip -\lastskip%
    \nobreak\vskip .5cm%
\fi
\else%
    \nobreak\vskip -\lastskip%
    \nobreak\vskip .5cm%
\fi%
\nextheadno%
\ifmaketdm%
    \def\temp{#1}%
    \edef\sectiontitle{\expandafter\gobbleeight\meaning\temp}%
    \edef\textoutput{\write\tdm{\noindent{\noexpand\mytdmheadfont\quad\headinfo\ - \sectiontitle\hfill\noexpand\folio}\par}}%
    \textoutput%
\fi%
\font\headfontA=cmbx10 at 14pt%
{\headfontA\noindent\headinfo\ -\ #1.\hfil}%
\nobreak\vskip .5cm%
}%


\def\newsubhead#1%
{%
\ifhmode%
    \mypar%
\fi%
\ifnum\subheadno=0%
\else%
    \penalty\headpenalty\vskip .4cm%
\fi%
\nextsubheadno%
\ifmaketdm%
    \def\temp{#1}%
    \edef\subsectiontitle{\expandafter\gobbleeight\meaning\temp}%
    \edef\textoutput{\write\tdm{\noindent{\noexpand\mytdmsubheadfont\quad\quad\subheadinfo\ - \subsectiontitle\hfill\noexpand\folio}\par}}%
    \textoutput%
\fi%
\font\subheadfontA=cmsl10 at 12pt
{\subheadfontA\noindent\subheadinfo\ #1.\hfil}%
\nobreak\vskip .25cm %
}%

%
%


\font\mathromanten=cmr10
\font\mathromanseven=cmr7
\font\mathromanfive=cmr5
\newfam\mathromanfam
\textfont\mathromanfam=\mathromanten
\scriptfont\mathromanfam=\mathromanseven
\scriptscriptfont\mathromanfam=\mathromanfive
\def\mathroman{\fam\mathromanfam}


\font\sansseriften=cmss10
\font\sansserifseven=cmss7
\font\sansseriffive=cmss5
\newfam\sansseriffam
\textfont\sansseriffam=\sansseriften
\scriptfont\sansseriffam=\sansserifseven
\scriptscriptfont\sansseriffam=\sansseriffive
\def\mathsf{\fam\sansseriffam}


\font\boldten=cmb10
\font\boldseven=cmb7
\font\boldfive=cmb5
\newfam\mathboldfam
\textfont\mathboldfam=\boldten
\scriptfont\mathboldfam=\boldseven
\scriptscriptfont\mathboldfam=\boldfive
\def\mathbf{\fam\mathboldfam}


\font\mycmmiten=cmmi10
\font\mycmmiseven=cmmi7
\font\mycmmifive=cmmi5
\newfam\mycmmifam
\textfont\mycmmifam=\mycmmiten
\scriptfont\mycmmifam=\mycmmiseven
\scriptscriptfont\mycmmifam=\mycmmifive

\def\hexa#1{\ifcase #1 0\or 1\or 2\or 3\or 4\or 5\or 6\or 7\or 8\or 9\or A\or B\or C\or D\or E\or F\fi}
\mathchardef\mathi="7\hexa\mycmmifam7B
\mathchardef\mathj="7\hexa\mycmmifam7C


\font\mymsbmten=msbm10 at 8pt
\font\mymsbmseven=msbm7 at 5.6pt
\font\mymsbmfive=msbm5 at 4pt
\newfam\mymsbmfam
\textfont\mymsbmfam=\mymsbmten
\scriptfont\mymsbmfam=\mymsbmseven
\scriptscriptfont\mymsbmfam=\mymsbmfive

\mathchardef\mybeth="7\hexa\mymsbmfam69
\mathchardef\mygimmel="7\hexa\mymsbmfam6A
\mathchardef\mydaleth="7\hexa\mymsbmfam6B


\def\placelabel[#1][#2]#3{{%
\setbox10=\hbox{\raise #2cm \hbox{\hskip #1cm #3}}%
\ht10=0pt%
\dp10=0pt%
\wd10=0pt%
\box10}}%


\newif\ifinproclaim%
\global\inproclaimfalse%
\def\proclaim#1{%
\medskip%
%
%
\bgroup%
\inproclaimtrue%
\setbox10=\vbox\bgroup\leftskip=0.8em\noindent{\bf #1}\sl%
}

\def\endproclaim{%
\egroup%
\setbox11=\vtop{\noindent\vrule height \ht10 depth \dp10 width 0.1em}%
\wd11=0pt%
\setbox12=\hbox{\copy11\kern 0.3em\copy11\kern 0.3em}%
\wd12=0pt%
\setbox13=\hbox{\noindent\box12\box10}%
\noindent\unhbox13%
\egroup%
\medskip\ignorespaces%
}

\def\proclaim#1{%
\medskip%
\bgroup%
\inproclaimtrue%
\noindent{\bf #1}%
\nobreak\medskip%
\sl%
}

\def\endproclaim{%
\mypar\egroup\penalty\proclaimpenalty\medskip\ignorespaces%
}

\def\noskipproclaim#1{%
\medskip%
\bgroup%
\inproclaimtrue%
\noindent{\bf #1}\nobreak\sl%
}

\def\endnoskipproclaim{%
\mypar\egroup\penalty\proclaimpenalty\medskip\ignorespaces%
}


\def\ninn{{n\in\Bbb{N}}}

\def\proof{{\noindent\bf Proof:\ }}

\def\remark{{\noindent\sl Remark:\ }}

\def\mlim{\mathop{{\mathroman Lim}}}

\def\msup{\mathop{{\mathroman Sup}}}
\def\minf{\mathop{{\mathroman Inf}}}
\def\msf#1{{\mathsf #1}}

\def\qed{~$\square$}
\def\munion{\mathop{\cup}}
\def\minter{\mathop{\cap}}
\def\myitem#1{%
    \noindent\hbox to .5cm{\hfill#1\hss}
}

\catcode`\@=11
\def\Eqalign#1{\null\,\vcenter{\openup\jot\m@th\ialign{%
\strut\hfil$\displaystyle{##}$&$\displaystyle{{}##}$\hfil%
&&\quad\strut\hfil$\displaystyle{##}$&$\displaystyle{{}##}$%
\hfil\crcr #1\crcr}}\,}
\catcode`\@=12

\def\makeop#1{%
\global\expandafter\def\csname op#1\endcsname{{\mathroman #1}}}%

\def\makeopsmall#1{%
\global\expandafter\def\csname op#1\endcsname{{\mathroman{\lowercase{#1}}}}}%

\makeopsmall{ArcTan}%
\makeopsmall{ArcCos}%
\makeop{Arg}%
\makeop{Det}%
\makeop{Log}%
\makeop{Re}%
\makeop{Im}%
\makeop{Dim}%
\makeopsmall{Tan}%
\makeop{Ker}%
\makeopsmall{Cos}%
\makeopsmall{Sin}%
\makeop{Exp}%
\makeopsmall{Tanh}%
\makeop{Tr}%
\makeop{Mob}
\makeop{End}%
\makeop{Long}%
\makeop{Ch}%
\makeop{Exp}%
\makeop{Ln}%
\makeop{PSO}
\makeop{PSL}
\makeop{Int}%
\makeop{Ext}%
\makeop{Aire}%
\makeop{Im}%
\makeop{Conf}%
\makeop{Exp}%
\makeop{Mod}%
\makeop{Log}%
\makeop{Ext}%
\makeop{Int}%
\makeop{Dist}%
\makeop{Aut}%
\makeop{Id}%
\makeop{SO}%
\makeop{min}
\makeopsmall{Coth}
\makeop{Homeo}%
\makeop{Vol}%
\makeop{Ric}%
\makeop{Hess}%
\makeop{Euc}%
\makeop{Isom}%
\makeop{Max}%
\makeop{Long}%
\makeop{Inj}
\makeop{Fixe}%
\makeop{Wind}%
\makeop{Imm}
\makeop{Acc}
\makeop{Mush}%
\makeop{Hom}
\makeop{PSO}
\makeop{Ad}%
\makeop{loc}%
\makeop{Len}%
\makeop{Pleat}
\makeop{Area}%
\makeop{Diam}
\makeop{Rep}
\makeop{SL}%
\makeop{GL}%
\makeop{dVol}%
\makeop{Min}%
\makeop{Symm}%
\makeop{Met}
\makeop{O}%
\makeop{Dev}
\makeopsmall{LimSup}
\makeopsmall{CoTanh}
\makeopsmall{ArcTanh}
\makeopsmall{Sinh}
\makeopsmall{Cosh}
\makeopsmall{ArcCoTanh}

\let\emph=\bf

\hyphenation{quasi-con-formal}

%
%

\ifmakebiblio%
    \openout\biblio=biblio.tex %
    {%
        \edef\fileoutput{\write\biblio{\bgroup\leftskip=2em}}%
        \fileoutput
    }%
\fi%

\newref{AndBarbBegZegh}{Andersson L., Barbot T., B\'e guin F., Zeghib A., Cosmological time versus CMC time in spacetimes of constant curvature}
\newref{Aubin}{Aubin T., {\sl Nonlinear analysis on manifolds. Monge-Amp\`ere equations}, Die Grund\-lehren der mathematischen Wissenschaften, {\bf 252}, Springer-Verlag, New York, (1982)}
\newref{CaffNirSpr}{Caffarelli L., Nirenberg L., Spruck J., The Dirichlet problem for nonlinear second-order elliptic equations. I. Monge-Amp\`ere equation. {\sl Comm. Pure Appl. Math} {\bf 37} (1984), no. 3, 369--402} 
\newref{ChenGreen}{Chen S., Greenberg L., Hyperbolic Spaces, Contribution to {\sl Analysis}, Academic Press, New York, (1974), 49-87}
\newref{EpsMard}{Epstein D. B. A., Marden, A., Convex hulls in hyperbolic space, a theorem of Sullivan, and measured pleated surfaces, 
In {\sl Fundamentals of hyperbolic geometry: selected expositions}, London Math. Soc. Lecture Note Ser., {\bf 328}, Cambridge Univ. Press, Cambridge, (2006)}
\newref{Fried}{Fried D., Closed Similarity Manifolds, {\sl Comment. Math. Helvetici} {\bf 55} (1980), 576-582}
\newref{GuanSpruck}{Guan B., Spruck J., The existence of hypersurfaces of constant Gauss curvature with prescribed boundary, {\sl J. Differential Geom.} {\bf 62} (2002), no. 2, 259--287}
\newref{HarveyLawson}{Harvey R., Lawson H. B. Jr., Calibrated geometries, {\sl Acta. Math.} {\bf 148} (1982), 47--157}
\newref{KamiA}{Kamishima T., Conformally Flat Manifolds whose Development Maps are not Surjective, {\sl Trans. Amer. Math. Soc.} {\bf 294} (1986), no. 2, 607-623}
\newref{KamTan}{Kamishima Y., Tan S., Deformation spaces on geometric structures, In {\sl Aspects of low-dimensional manifolds}, Adv. Stud. Pure Math., {\bf 20}, Kinokuniya, Tokyo, (1992)}
\newref{KapB}{Kapovich M., Deformation spaces of flat conformal structures. Proceedings of the Second Soviet-Japan Joint Symposium of 
Topology (Khabarovsk, 1989), {\sl Questions Answers Gen. Topology} {\bf 8} (1990), no. 1, 253--264}
\newref{KraSch}{Krasnov K., Schlenker J.M., On the renormalized volume of hyperbolic $3$-manifolds, math.DG/0607081} 
\newref{Kulkarni}{Kulkarni R.S., Pinkall U., A canonical metric for M\"obius structures and its applications, {\sl Math. Z.} {\bf 216} (1994), no.1, 89--129}
\newref{LabA}{Labourie F., Un lemme de Morse pour les surfaces convexes, {\sl Invent. Math.} {\bf 141} (2000), 239--297}
\newref{LabB}{Probl\`eme de Minkowski et surfaces \`a courbure constante dans les vari\'et\'es hyperboliques, {\sl Bull. Soc. Math. Fr.} {\bf 119} (1991), 307-325}
\newref{MazzPac}{Mazzeo R., Pacard F., Constant curvature foliations in asymptotically hyperbolic spaces}
\newref{SmiA}{Smith G., Special Legendrian structures and Weingarten problems, Preprint, Orsay (2005)}
\newref{SmiNLD}{Smith G., The non-linear Dirichlet Problem in Hadamard Manifolds, in preparation}
\newref{SmiH}{Smith G., A brief note on foliations of constant Gaussian curvature, arXiv:0802.2202}
\newref{Thurston}{Thurston W., {\sl Three-dimensional geometry and topology}, Princeton Mathematical Series, {\bf 35}, Princeton University Press, Princeton, NJ, (1997)}
%
\ifmakebiblio%
    {\edef\fileoutput{\write\biblio{\egroup}}%
    \fileoutput}%
\fi%

%
%
%
\document
\myfontdefault
\global\chapno=1
\global\showpagenumflag=1
\def\Pagetitle{}
\null
\vfill
\def\centre{\rightskip=0pt plus 1fil \leftskip=0pt plus 1fil \spaceskip=.3333em \xspaceskip=.5em \parfillskip=0em \parindent=0em}%
\def\textmonth#1{\ifcase#1\or January\or Febuary\or March\or April\or May\or June\or July\or August\or September\or October\or November\or December\fi}
\font\abstracttitlefont=cmr10 at 14pt
{\abstracttitlefont\centre Moduli of Flat Conformal Structures of Hyperbolic Type\par}
\bigskip
{\centre Graham Smith\par}
\bigskip
{\centre \the\day\ \textmonth\month\ \the\year\par}
\bigskip
{\centre Centre de Recerca Matem\`atica,\par
Facultat de Ci\`encies, Edifici C,\par
Universitat Aut\`onoma de Barcelona,\par
08193 Bellaterra,\par
Barcelona,\par
SPAIN\par}
\bigskip
\noindent{\emph Abstract:\ }To each flat conformal structure (FCS) of hyperbolic type in the sense of Kulkarni-Pinkall, we associate, for all $\theta\in[(n-1)\pi/2,n\pi/2[$ and for all $r>\opTan(\theta/n)$ a unique immersed hypersurface $\Sigma_{r,\theta}=(M,i_{r,\theta})$ in $\Bbb{H}^{n+1}$ of constant $\theta$-special Lagrangian
curvature equal to $r$. We show that these hypersurfaces smoothly approximate the boundary of the canonical hyperbolic end associated to the FCS by Kulkarni and Pinkall and thus obtain results concerning the continuous dependance of the hyperbolic end and of the Kulkarni-Pinkall metric on the flat conformal structure.
\par 
\par 
\vfill
\nextoddpage
%
%
\def\Pagetitle{\sl Moduli of Flat Conformal Structures of Hyperbolic Type}
\global\pageno=1
\newhead{Introduction}
\noindent A flat conformal structure (FCS) (or M\"obius structure) on an $n$-dimensional manifold, $M$, is an atlas of $M$ whose charts lie in $S^n$ and whose transition maps are restrictions of conformal (i.e. M\"obius) mappings of $S^n$. Such structures arise naturally in different domains of mathematics. To every FCS of hyperbolic type may be canonically associated a complete hyperbolic manifold with convex boundary called the hyperbolic end of that structure. The purpose of this paper is to associate to every such FCS defined over a compact manifold families of foliations of neighbourhoods of the finite boundary of its hyperbolic end consisting of smooth, convex hypersurfaces of constant curvature.
\medskip
\noindent The history of FCSs begins with the $2$-dimensional case. Here, Thurston shows, for example, that the moduli space of FCSs over a compact surface, $M$, is homeomorphic to the Cartesian product $\Cal{T}\times\Cal{ML}(M)$ of the Teichm\"uller space of $M$ with the space of measured geodesic laminations over $M$ (see \cite{KamTan} or \cite{Thurston} for details). An important step in Thurston's proof involves the construction of a convex, pleated, equivariant 
``immersion'' $i_T:\tilde{M}\rightarrow\Bbb{H}^3$ from the universal cover of $M$ into $\Bbb{H}^3$ which is canonically associated to the FCS. This construction generalises that of the Nielsen Kernel of a quasi-Fuchsian manifold (see \cite{EpsMard} for a detailed study of its properties in this case).
\medskip
\noindent In the higher dimensional case, Kapovich \cite{KapB} provides information on the moduli space of FCSs, but much remains unknown. However, when $M$ is of hyperbolic type (see section \headref{ConformalStructures}), Kulkarni and Pinkall showed in \cite{Kulkarni} that Thurston's construction may still be carried out. This yields a convex,
stratified, equivariant ``immersion'' $i_{KP}:M\rightarrow\Bbb{H}^{n+1}$ in $\Bbb{H}^{n+1}$ canonically associated to the M\"obius structure, as well as a canonical $C^{1,1}$ metric over $M$ with a.e. defined sectional curvatures taking values between $-1$ and $1$. We call this metric the Kulkarni-Pinkall metric of the M\"obius structure and denote it by $g_{KP}$.
\medskip
\noindent Heuristically, a hyperbolic end over a manifold $M$ is a complete, hyperbolic manifold with concave, stratified boundary whose interior is homeomorphic to $M\times\Bbb{R}$. A detailed description is provided in Sections \headref{HyperbolicEnds} and \headref{HeadGBC}. Strictly speaking, we call the boundary of $\Cal{E}$ the finite boundary, and we denote it by $\partial_0\Cal{E}$. This distinguishes it from the ideal boundary, $\partial_\infty\Cal{E}$, which is the set of equivalence classes of complete half geodesics whose distance from $\partial_0\Cal{E}$ tends to infinity.
\medskip
\noindent In \cite{Kulkarni}, Kulkarni and Pinkall show that the ``immersion'' $i_{KP}$ may be interpreted as the finite boundary of a hyperbolic end, $\Cal{E}$ which is also canonically associated to the FCS and whose 
ideal boundary $\partial_\infty\Cal{E}$ is conformally equivalent to $M$. $\Cal{E}$ thus provides a cobordism between
$i_{KP}$ and $M$. It is for neighbourhoods of the finite boundaries of these hyperbolic ends that we construct foliations by hypersurfaces of constant curvature. These foliations may thus be considered as families of smoothings of $i_{KP}$. This construction generalises to higher dimensions the result \cite{LabB} of Labourie which provides families of parametrisations of the moduli spaces of three dimensional hyperbolic manifolds with geometrically finite ends.
\medskip
\noindent The special Lagrangian curvature, $R_\theta$ was first developed by the author in \cite{SmiA}. We recall its 
construction in section \headref{SpecialLagrangianCurvature}. Its most important properties are that it is only defined for strictly convex immersed hypersurfaces and that it is regular in a PDE sense, which is summarised in this paper in terms of Theorems \procref{CompactnessA} and \procref{CompactnessB} (proven in \cite{SmiA}) and Theorem \procref{TheoremLocalExistence} (proven in \cite{SmiNLD}).
\medskip
\noindent Of tangential interest, this notion of curvature arises from the natural special Legendrian structure of the unitary bundle of $U\Bbb{H}^3$. Special Legendrian structures are closely related to special Lagrangian structures which are studied under the heading of Calabi-Yau manifolds. Special Lagrangian and Legendrian submanifolds have themselves been of growing interest to mathematicians and physicists since the landmark paper \cite{HarveyLawson} of Harvey and Lawson concerning calibrated geometries. In its classical form, the special Lagrangian operator is a second order,
highly non-linear partial differential operator of determinant type closely related to the Monge-Amp\`ere operator, and which is among the archetypical highly non-linear partial differential operators studied in detail in most standard works on nonlinear PDEs (\cite{Aubin} and \cite{CaffNirSpr} to name but two).
\medskip
\noindent The main results of this paper are most appropriately described in terms of developing maps (see section \headref{ConformalStructures}). Let $M$ be a manifold. A M\"obius structure over $M$ may be considered as a pair $(\varphi,\theta)$ where $\theta:\pi_1(M)\rightarrow\opConf(S^n)$ is a homomorphism and $\varphi:\tilde{M}\rightarrow S^n$ is a local homeomorphism from the universal cover of $M$ into $S^n$ which is equivariant with respect to $\theta$. Two pairs are equivalent if and only if they differ by a conformal mapping of $S^n$. We furnish the space of M\"obius structures with the (quotient of) the topology of local uniform convergence. $\varphi$ is called the developing map and $\theta$ is called the holonomy of the M\"obius structure.
\medskip
\noindent We define the Gauss mapping $\overrightarrow{n}:U\Bbb{H}^{n+1}\rightarrow\partial_\infty\Bbb{H}^{n+1}$ as follows. For $v$ a unit vector in $U\Bbb{H}^{n+1}$, let $\gamma_v:[0,+\infty[\rightarrow\Bbb{H}^{n+1}$ be the half geodesic such that $\partial_t\gamma(0)=v$. We define:
$$
\overrightarrow{n}(v) = \gamma_v(+\infty) = \mlim_{t\rightarrow+\infty}\gamma_v(+\infty).
$$
\noindent Let $i:M\rightarrow\Bbb{H}^{n+1}$ be a convex immersion. Since $i$ is convex, there exists a unique exterior vector field $\msf{N}_i$ over $i$ in $U\Bbb{H}^{n+1}$. We say that $i$ {\bf projects asymptotically} to the M\"obius structure $(\varphi,\theta)$ if and only if $i$ is equivariant with respect to $\theta$, and, up to reparametrisation:
$$
\overrightarrow{n}\circ \msf{N}_i = \varphi.
$$
\proclaim{Theorem \nextprocno}
\noindent Let $\eta\in](n-1)\pi/2,n\pi/2[$ be an angle, and let $r>\opTan(\eta/n)$. Let $M$ be a compact $n$ dimensional manifold and let $(\varphi,\theta)$ be an FCS of hyperbolic type over $M$. There exists a unique, convex, equivariant immersion $i_{r,\eta}:\tilde{M}\rightarrow\Bbb{H}^{n+1}$ such that:
\medskip
\myitem{(i)} $i_{r,\eta}$ is a graph over $i_{KP}$;
\medskip
\myitem{(ii)} $i_{r,\eta}$ projects asymptotically to $\varphi$; 
\medskip
\myitem{(iii)} $R_\eta(i_{r,\eta}) = r$.
\medskip
\noindent Moreover, if $(\varphi,\theta)$ is not conformally equivalent to $S^{n-1}\times S^1$, where $S^k$ is the $k$-dimensional sphere, then the same result holds for $\eta=(n-1)\pi/2$.
\endproclaim
\proclabel{TheoremExistenceAndUniqueness}
\remark The proof of this theorem uses the Perron method. The finite boundary forms a barrier, which follows from the Geodesic Boundary Property (see Definition \procref{DefinitionGBC}). In particular, as in the remarks following Definition \procref{DefinitionGBC}, the existence result in fact holds in a much more general class of negatively curved ends of non-constant sectional curvature bounded above by $-1$ whose finite boundary possesses the Geodesic Boundary Property or even the weak Geodesic Boundary Property.
\medskip
\noindent Since they are graphs over the Kulkarni-Pinkall immersion, these immersed hypersurfaces may be considered as submanifolds of the hyperbolic end of the FCS:
\proclaim{Theorem \nextprocno}
\noindent Let $\Cal{E}$ be the hyperbolic end of an FCS. Let $\theta\in[(n-1)\pi/2,n\pi/2[$ be an angle. For all $r>\opTan(\theta/n)$, let $\Sigma_{r,\theta}=(S,i_{r,\theta})$ be the unique, smooth, convex, immersed hypersurface on $\Cal{E}$ which is a graph over $\partial\Cal{E}$ and which satisfies $R_\theta(i_{r,\theta})=r$.
\medskip
\noindent The family $(\Sigma_{r,\theta})_{r>\opTan(\theta/n)}$ foliates a neighbourhood, $\Omega_\theta$, of $\partial\Cal{E}$. Morever $(\hat{\Sigma}_{r,\theta})_{r>\opTan(\theta/n)}$ converges towards $N\Cal{E}$ in the $C^0$ sense as $r$ tends to $+\infty$, and, for any compact 
subset, $K$, of $\Cal{E}$, there exists $\theta_0<n\pi/2$ such that for $\theta>\theta_0$, $K\subseteq\Omega_\theta$.
\endproclaim
\proclabel{TheoremFoliations}
\remark The final part of this theorem suggests that by judiciously choosing $r$ as a function of $\theta$, it may be possible to obtain smooth foliations of the entire hyperbolic end.
\medskip
\remark Towards completion of this paper, the author was made aware of a recent, complementary result of Mazzeo and Pacard \cite{MazzPac}. There, using entirely different techniques, and under different hypotheses on the hyperbolic end, the authors prove the existence of foliations by constant mean curvature hypersurfaces near the ideal boundary, though not near the finite boundary, as is obtained here. It appears reasonable that a happy marriage of these techniques could yield more detailed information concerning the structure of the hyperbolic end and its relation to its ideal boundary.
\medskip
\noindent In the special case where $\Cal{E}$ is an end of a quasi-Fuchsian manifold, the foliations may be extended up to the ideal boundary, and we obtain:
\proclaim{Theorem \nextprocno}
\noindent Let $\Cal{E}$ be a hyperbolic end of a quasi-Fuchsian manifold. Let $\theta\in[(n-1)\pi/2,n\pi/2[$ be an angle. For all $r>\opTan(\theta/n)$, let $\Sigma_{r,\theta}=(S,i_{r,\theta})$ be the unique, smooth, convex, immersed hypersurface on $\Cal{E}$ which is a graph over $\partial\Cal{E}$ and which satisfies $R_\theta(i_{r,\theta})=r$.
\medskip
\noindent The family $(\Sigma_{r,\theta})_{r>\opTan(\theta/n)}$ foliates $\Cal{E}$. Morever $(\hat{\Sigma}_{r,\theta})_{r>\opTan(\theta/n)}$ converges towards $N\Cal{E}$ in the $C^0$ sense as $r$ tends to $+\infty$, and $(\Sigma_{r,\theta})_{r>\opTan(\theta/n)}$ converges to $\partial_\infty\Cal{E}$ in the Hausdorff sense as $r$ tends to $\opTan(\theta/n)$.
\endproclaim
\proclabel{TheoremFoliationsInQuasiFuchsian}
\remark In fact, this result holds for any FCS whose developing map avoids an open subset of $\partial_\infty\Bbb{H}^{n+1}$.
\medskip
\noindent We next consider how these foliations vary with the FCS:
\goodbreak
\proclaim{Theorem \nextprocno}
\noindent Let $M$ be a compact manifold. Let $(\theta_x,\varphi_x)_{\|x\|<\epsilon}$ be a continuous family of FCSs of hyperbolic type over $M$ whose holonomy varies smoothly. Let $\theta\in[(n-1)\pi/2,n\pi/2[$ be an angle, and let $r>\opTan(\theta/n)$. For all $x$, let $\Sigma_x=(S,i_x)$ be the unique, smooth, convex, immersed hypersurface in $\Cal{E}(\theta_x,\varphi_x)$ such that $R_\theta(i_x)=r$. Then, up to reparametrisation, $i_x$ varies smoothly with $x$.
\endproclaim
\proclabel{TheoremContinuousDependance}
\remark It follows that the space of hypersurfaces of constant special Lagrangian curvature yields smooth moduli for the
space of FCSs of hyperbolic type over $M$ which are compatible with the smooth structure obtained from the canonical embedding of this space into $\opPSO(n+1,1)^{\pi_1(M)}$, and which also, importantly, encode smooth information about the hyperbolic end and the Kulkarni-Pinkall metric.
\medskip
\noindent As an illustration of these results, we now consider two special cases. The first is when $n$ is equal to $2$, and $\theta=\pi/2$. Here the special Lagrangian curvature reduces to the Gaussian curvature and we recover the following,
now classical, result of Labourie \cite{LabB}:
\proclaim{Theorem \nextprocno, {\bf Labourie (1991)}}
\noindent Let $\Sigma$ be a compact surface of hyperbolic type. Let $(\alpha,\varphi)$ be an FCS over $\Sigma$ and let $\Cal{E}$ be the hyperbolic end of $(\alpha,\varphi)$. There exists a unique, smooth foliation $(\Sigma_k)_{k\in]0,1[}$ of $\Cal{E}$ such that:
\medskip
\myitem{(i)} for each $k$, $\Sigma_k$ is a smooth, immersed surface of constant Gaussian (extrinsic) curvature equal to $k$;
\medskip
\myitem{(ii)} $\Sigma_k$ tends to $\partial_0\Cal{E}$ in the Hausdorff sense as $k$ tends to $0$; and
\medskip
\myitem{(iii)} $\Sigma_k$ tends to $\partial_\infty\Cal{E}$ in the Hausdorff sense as $k$ tends to $1$.
\endproclaim
\remark The geometric properties particular to this special case allow us to extend the foliations up to the ideal boundary (see also \cite{MazzPac} and \cite{SmiH}).
\medskip
\noindent The second special case is when $n=3$ and $\theta=\pi$. In this case, the special Lagrangian curvature still has a very simple expression:
\proclaim{Theorem \nextprocno}
\noindent Let $M$ be a compact three dimensional manifold. Let $(\alpha,\varphi)$ be an FCS over $M$ of hyperbolic type. Let $\Cal{E}$ be the hyperbolic end of $(\alpha,\varphi)$. There exists a unique, smooth foliation 
$(\Sigma_r)_{r\in]3,+\infty[}$ of $\Cal{E}$ such that:
\medskip
\myitem{(i)} for each $r$, $\Sigma_r$ is a smooth, immersed hypersurface such that:
$$
H(\Sigma_r)/K(\Sigma_r) = r,
$$
\noindent where $H(\Sigma_r)$ and $K(\Sigma_r)$ are the mean and Gaussian curvatures of $\Sigma_r$ respectively; and
\medskip
\myitem{(ii)} $\Sigma_r$ tends to $\partial_0\Cal{E}$ in the Hausdorff sense as $r$ tends to $+\infty$.
\endproclaim
\proclabel{TheoremFourDSpecialCase}
\noindent Towards completion of this paper, the author was made aware of related work by Andersson, Barbot, B\'eguin and Zeghib \cite{AndBarbBegZegh}. Here the authors study constant mean curvature foliations of Lorentzian, anti de-Sitter and de-Sitter spacetimes. There is a natural duality between hyperbolic ends and de-Sitter spacetimes, and thus a duality between their framework and our own. One interesting consequence is that, in the $4$-dimensional case, Theorem \procref{TheoremFourDSpecialCase} yields foliations of neighbourhoods of the past ends of four dimensional de-Sitter spacetimes by $3$-dimensional space-like hypersurfaces of constant scalar curvature. This may be related to the Yamabe problem of the flat conformal structure, which is relevant to \cite{MazzPac}.
\medskip
\noindent Finally, the proofs of these theorems requires a detailed understanding of the geometric structure of the Kulkarni-Pinkall hyperbolic end of a flat conformal structure. We obtain the following characterisation of the Kulkarni-Pinkall end in terms of completeness and local geometric data, which the author is not aware of in the litterature:
\proclaim{Theorem \nextprocno}
\noindent Let $\tilde{N}$ be a hyperbolic end. Suppose that:
\medskip
\myitem{(i)} $\tilde{N}$ possesses the Geodesic Boundary Property; and
\medskip
\myitem{(ii)} $\tilde{N}$ is complete.
\medskip
\noindent Then $\tilde{N}$ is the Kulkarni-Pinkall hyperbolic end of its quotient M\"obius manifold.
\medskip
\noindent Moreover, if $N$ is a compact M\"obius manifold, then the family of hyperbolic ends whose quotient M\"obius manifold is $N$ is partially ordered, and the Kulkarni-Pinkall hyperbolic end of $N$ is the unique maximal element of this family. 
\endproclaim
\proclabel{TheoremStructureOfKPEnd}
\noindent Indeed, as noted in the remark following Theorem \procref{TheoremContinuousDependance}, the foliations constructed here encode smooth information about the hyperbolic end whilst depending smoothly on the conformal structure. We therefore expect them to be of considerable use in the future study of FCSs. Indeed, as examples of possible applications of these results, we state two immediate corollaries. The first concerns continuous dependence of $i_{KP}$:
\proclaim{Theorem \nextprocno}
\noindent Let $M$ be a compact manifold. Let $(\theta_n,\varphi_n)_\ninn,(\theta_0,\varphi_0)$ be FCSs of hyperbolic type over $M$ such that $(\theta_n,\varphi_n)_\ninn$ converges to $(\theta_0,\varphi_0)$, then
$(N\partial_0\Cal{E}(\theta_n,\varphi_n))_\ninn$ converges to $(N\partial_0\Cal{E}(\theta_0,\varphi_0))$ in the $C^0$ sense.
\endproclaim
\noindent And the second result concerns the Kulkarni-Pinkall metric. Let $D$, $V$ and $I$ represent the diameter, volume and injectivity radius respectively of the Kulkarni-Pinkall metric. We obtain the following continuity and compactness result:
\proclaim{Theorem \nextprocno}
\noindent Let $M$ be a compact manifold. $D$, $V$ and $I$ define continuous functions over the space of FCSs of hyperbolic type over $M$. Moreover, the pairs $(I,D)$ and $(I,V)$ define proper functions over the space of FCSs of hyperbolic type.
\endproclaim
\noindent This paper is structured as follows:
\medskip
\myitem{(a)} In Sections $2$ to $6$, we introduce the various concepts used in this paper. In particular, hyperbolic ends and their relationship to flat conformal structures are studied in sections $4$ to $6$ and Theorem \procref{TheoremStructureOfKPEnd} is proven in section $6$;
\medskip
\myitem{(b)} In Section $7$ \& $8$, we describe how equivariant immersions of smooth curvature may be deformed, which is relevant to Theorem \procref{TheoremContinuousDependance};
\medskip
\myitem{(c)} In Section $9$, we recall the compactness properties of the special Lagrangian curvature, which form an important component of the Perron method used to prove Theorem \procref{TheoremExistenceAndUniqueness};
\medskip
\myitem{(d)} In Section $10$, we recall the geometric maximum principal which is used to control the location of immersed hypersurfaces of given SL curvature;
\medskip
\myitem{(e)} In Section $11$, we prove the uniqueness part of Theorem \procref{TheoremExistenceAndUniqueness};
\medskip
\myitem{(f)} Theorems \procref{TheoremExistenceAndUniqueness}, \procref{TheoremFoliations} and \procref{TheoremContinuousDependance} are proven in Section $12$;
\medskip
\myitem{(h)} In Section $13$, quasi-Fuchsian manifolds are introduced and Theorem \procref{TheoremFoliationsInQuasiFuchsian} is proven; and
\medskip
\myitem{(i)} In Appendix $A$, we show how the Kulkarni-Pinkall metric may be used to furnish a simpler proof of a result of Kamishima.
\medskip
\noindent This paper has known a long and tortuous evolution since its conception. I would like to thank Kirill Krasnov, Fran\c{c}ois Labourie and Jean-Marc Schlenker for encouraging me to study this problem in the first place. I am equally grateful to Werner Ballmann, Ursula Hamenstaedt and Joan Porti for many useful conversations about FCSs (and to the latter two for drawing attention to the various errors in earlier drafts of this paper). Finally, I would like to thank the Max Planck Institutes for Mathematics in the Sciences in Leipzig, the Max Planck Institute for Mathematics in Bonn and the Centre de Recerca Matem\`atica in Barcelona for providing the conditions required to carry out this research.
\goodbreak
\newhead{Immersed Submanifolds and the Cheeger/Gromov Topology}
\noindent Let $M$ be a smooth Riemannian manifold. An {\bf immersed submanifold} is a pair 
$\Sigma=(S,i)$ where $S$ is a smooth manifold and $i:S\rightarrow M$ is a smooth immersion. A 
{\bf pointed immersed submanifold} in  $M$ is a pair $(\Sigma,p)$ where $\Sigma=(S,i)$ is an immersed submanifold in $M$ and $p$ is a point in $S$. An {\bf immersed
hypersurface} is an immersed submanifold
of codimension $1$. We give $S$ the unique Riemannian metric $i^*g$ which makes $i$ into an isometry. We say that $\Sigma$ is {\bf complete} if and only if the Riemannian manifold $(S,i^*g)$ is.
\medskip
\noindent Let $UM$ be the unitary bundle of $M$ (i.e the bundle of unit vectors in $TM$. In the cooriented case (for example, when $I$ is convex), there exists a unique exterior normal vector 
field $\msf{N}$ over $i$. We denote $\hat{\mathi}=\msf{N}$ and call it the {\bf Gauss lift} of $i$. Likewise, we call the manifold $\hat{\Sigma}=(S,\hat{\mathi})$ the
{\bf Gauss lift} of $\Sigma$.
\medskip
\noindent A pointed Riemannian manifold is a pair $(M,p)$ where $M$ is a Riemannian manifold and $p$ is a point in $M$. Let $(M_n,p_n)_{n\in\Bbb{N}}$ be a sequence of pointed Riemannian manifolds. For all $n$, we denote by $g_n$ the Riemannian metric over $M_n$. We say that the
sequence $(M_n,p_n)_{n\in\Bbb{N}}$ converges to the pointed manifold $(M_0,p_0)$ in the {\bf Cheeger/Gromov} sense if and only if for all $n$, there exists a mapping 
$\varphi_n:(M_0,p_0)\rightarrow (M_n,p_n)$, such that, for every compact subset $K$ of $M_0$, there 
exists $N\in\Bbb{N}$ such that for all $n\geqslant N$:
\medskip
\myitem{(i)} the restriction of $\varphi_n$ to $K$ is a $C^\infty$ diffeomorphism onto its image, and
\medskip
\myitem{(ii)} if we denote by $g_0$ the Riemannian metric over $M_0$, then the sequence of metrics $(\varphi_n^*g_n)_{n\geqslant N}$ converges to $g_0$ in the
$C^\infty$ topology over $K$.
\medskip
\noindent We refer to the sequence $(\varphi_n)_{n\in\Bbb{N}}$ as a sequence of {\bf convergence mappings} of the sequence $(M_n, p_n)_{n\in\Bbb{N}}$ with respect to the limit $(M_0,p_0)$. The convergence
mappings are trivially not unique. 
\medskip
\noindent Let $(\Sigma_n,p_n)_{n\in\Bbb{N}}=(S_n,p_n,i_n)_{n\in\Bbb{N}}$ be a sequence of pointed
immersed submanifolds in $M$. We say that $(\Sigma_n,p_n)_{n\in\Bbb{N}}$ converges to 
$(\Sigma_0,p_0)=(S_0,p_0,i_0)$ in the {\bf Cheeger/Gromov sense} if and only if the sequence 
$(S_n,p_n)_{n\in\Bbb{N}}$ of underlying manifolds converges to $(S_0,p_0)$ in the Cheeger/Gromov sense, and, for every sequence
$(\varphi_n)_{n\in\Bbb{N}}$ of convergence mappings of $(S_n,p_n)_\ninn$ with respect to this limit, and
for every compact subset $K$ of $S_0$, the sequence of functions $(i_n\circ\varphi_n)_{n\geqslant N}$ converges to the function $(i_0\circ\varphi_0)$ in the $C^\infty$ topology
over $K$.
\goodbreak
\newhead{Special Lagrangian Curvature}
\noindent The special Lagrangian curvature, which only has meaning for strictly convex immersed hypersurfaces, is defined as follows. Denote by $\opSymm(\Bbb{R}^n)$ the space of symmetric matrices over $\Bbb{R}^n$. We define $\Phi:\opSymm(\Bbb{R}^n)\rightarrow\Bbb{C}^*$ by:
\headlabel{SpecialLagrangianCurvature}
$$
\Phi(A) = \opDet(I+iA).
$$
\noindent Since $\Phi$ never vanishes and $\opSymm(\Bbb{R}^n)$ is simply connected, there exists a unique analytic function $\tilde{\Phi}:\opSymm(\Bbb{R}^n)\rightarrow\Bbb{C}$
such that:
$$
\tilde{\Phi}(I) = 0,\qquad e^{\tilde{\Phi}(A)} = \Phi(A)\qquad\forall A\in\opSymm(\Bbb{R}^n).
$$
\noindent We define the function $\opArcTan:\opSymm(\Bbb{R}^n)\rightarrow(-n\pi/2,n\pi/2)$ by:
$$
\arctan(A) = \opIm(\tilde{\Phi}(A)).
$$
\noindent This function is trivially invariant under the action of $O(\Bbb{R}^n)$. If $\lambda_1, ..., \lambda_n$ are the eigenvalues of $A$, then:
$$
\opArcTan(A) = \sum_{i=1}^n\opArcTan(\lambda_i).
$$
\noindent For $r>0$, we define:
$$
\opSL_r(A) = \opArcTan(rA).
$$
\noindent If $A$ is positive definite, then $SL_r$ is a strictly increasing function of $r$. Moreover, $SL_0=0$ and $SL_\infty=n\pi/2$. Thus, for all $\theta\in]0,n\pi/2[$, there exists a unique $r>0$ such that:
$$
SL_r(A) = \theta.
$$
\noindent We define $R_\theta(A) = r$. $R_\theta$ is also invariant under the action of $O(n)$ on the space of positive definite, symmetric matrices.
\medskip
\noindent Let $M$ be an oriented Riemannian manifold of dimension $n+1$. Let $\Sigma=(S,i)$ be a strictly convex, immersed
hypersurface in $M$. For $\theta\in]0,n\pi/2[$, we define $R_\theta(\Sigma)$ (the {\emph $\theta$-special Lagrangian curvature} of $\Sigma$) by:
$$
R_\theta(\Sigma) = R_\theta(A_\Sigma),
$$
\noindent where $A_\Sigma$ is the shape operator of $\Sigma$.
\goodbreak
\newhead{Hyperbolic Ends}
\noindent For all $m$, let $\Bbb{H}^{m+1}$ be $(m+1)$-dimensional hyperbolic space. Let $U\Bbb{H}^{m+1}$ be the unitary bundle over 
$\Bbb{H}^{m+1}$. Let $K$ be a convex subset of $\Bbb{H}^{m+1}$. We define $\Cal{N}(K)$, the set of normals over $K$ by:
\headlabel{HyperbolicEnds}
$$
\Cal{N}(K) = \left\{ v_x\in U\Bbb{H}^{m+1}\text{ s.t. }x\in\partial K\text{ and $v_x$ is a supporting normal to $K$ at $x$.}\right\}
$$
\noindent $\Cal{N}(K)$ is a $C^{0,1}$ submanifold of $U\Bbb{H}^{m+1}$. Let $\Omega$ be an open subset of $\Cal{N}(K)$. We define $\Cal{E}(\Omega)$, the end over $\Omega$ by:
$$
\Cal{E}(\Omega) = \left\{\opExp(t v_x)\text{ s.t. }t\geqslant 0, v_x\in\Omega\right\}.
$$
\noindent We say that a subset of $\Bbb{H}^{m+1}$ has {\bf concave boundary} if and only if it is the end of some open subset of the set of normals of a convex set. We refer to $\Omega$ as the finite boundary of $\Cal{E}(\Omega)$.
\medskip
\noindent We extend this concept to more general manifolds. Let $(M,\partial M)$ be a smooth manifold with continuous boundary. A {\bf hyperbolic
end} 
over $M$ is an atlas $\Cal{A}$ such that:
\medskip
\myitem{(i)} every chart of $\Cal{A}$ has convex boundary, and
\medskip
\myitem{(ii)} the transition maps of $\Cal{A}$ are isometries of $\Bbb{H}^{m+1}$.
\medskip
\noindent We refer to $\partial M$ as the finite boundary of $M$. In the sequel, we will denote it by $\partial_0 M$ in order to differentiate it from the ideal boundary $\partial_\infty M$ of $M$.
\medskip
\noindent We can construct hyperbolic ends using continuous maps into $U\Bbb{H}^{m+1}$. Let $M$ be an $m$-dimensional manifold without boundary. Let 
$i:M\rightarrow U\Bbb{H}^{m+1}$ be a continuous map. We say that $i$ is a {\bf convex immersion} if and only if for every $p$ in 
$M$, there exists a neighbourhood $\Omega$ of $p$ in $M$ and a convex subset $K\subseteq\Bbb{H}^{m+1}$ such that the restriction of $i$ to $\Omega$ is a homeomorphism 
onto an open subset of $\Cal{N}(K)$. In this case, we define the mapping $I:M\times[0,\infty[\rightarrow\Bbb{H}^{m+1}$ by:
$$
I(p,t) = \opExp(t i(p)).
$$
\noindent We refer to $I$ as the end of $i$. $I$ is a local homeomorphism from $M\times]0,\infty[$ into $\Bbb{H}^{m+1}$. If $g$ is the hyperbolic metric over $\Bbb{H}^{m+1}$, then
$I^*g$ defines a hyperbolic metric over this interior. $I^*g$ degenerates over the boundary, and we identify points that may be joined by curves of zero length. We denote
this equivalence by $\sim$ and we define $\Cal{E}(i)$, which we also call the end of $i$ by:
$$
\Cal{E}(i) = (M\times]0,\infty[)\munion(M/\sim).
$$
\noindent We shall see presently that every hyperbolic end may be constructed in this manner. Thus, if $\hat{M}$ is an end, and if $i:M\rightarrow U\Bbb{H}^{m+1}$ is a convex immersion such that $\hat{M}=\Cal{E}(i)$, then we say that $i$ is the boundary immersion of $\hat{M}$.
\goodbreak
\newhead{Flat Conformal Structures}
\noindent Let $\Bbb{H}^{n+1}$ be $(n+1)$-dimensional hyperbolic space. We identify $\partial_\infty\Bbb{H}^{n+1}$ with the $n$-dimensional sphere
$S^n$. $\opIsom(\Bbb{H}^{n+1})$ is identified with $\opPSO(n+1,1)$. This group acts faithfully on $S^n=\partial_\infty\Bbb{H}^{n+1}$. The image is a subgroup
of the group of homeomorphisms of the sphere. We denote this group by $\opMob(n)$ and we call elements of $\opMob(n)$ {\bf conformal maps}.
\headlabel{ConformalStructures}
\medskip
\noindent Let $M$ be a manifold. A flat conformal structure (FCS) on $M$ is an atlas $\Cal{A}$ of $M$ in $S^n$ whose transformation maps are restrictions 
of elements of $\opMob(n)$. Trivially, every element of $\opMob(n)$ is uniquely determined by its germ at a point. Thus, any chart of $\Cal{A}$ 
uniquely extends to a local homeomorphism from $\tilde{M}$, the universal cover of $M$, into $S^n$ which is equivariant with respect to a given
homomorphism. This yields an alternative definition of FCSs which is better adapted to our purposes:
\proclaim{Definition \nextprocno}
\noindent Let $M$ be a manifold. Let $\pi_1(M)$ be its fundamental group and let $\tilde{M}$ be its universal cover. A {\bf flat conformal structure}
over $M$ is a pair $(\varphi,\theta)$ where:
\medskip
\myitem{(i)} $\theta:\pi_1(M)\rightarrow\opMob(n)$ is a homomorphism, and
\medskip
\myitem{(ii)} $\varphi:\tilde{M}\rightarrow S^n$ is a local homeomorphism which is equivariant with respect to $\theta$.
\medskip
\noindent $\theta$ is called the {\bf holonomy} and $\varphi$ is called the {\bf developing map} of the flat conformal structure.
\medskip
\noindent We refer to a pair $(M,(\varphi,\theta))$ consisting of a manifold $M$ and a flat conformal structure over $M$ as a {\bf M\"obius manifold}.
In the sequel, where no ambiguity arises, we refer to the manifold with its conformal structure merely by $M$.
\endproclaim
\remark A canonical differential structure on $M$ is obtained by pulling back the differential structure of $S^n$ through $\varphi$.
\medskip
\noindent M\"obius manifolds are divided into three types (for more details, see \cite{Kulkarni}):
\medskip
\myitem{(i)} manifolds of {\bf elliptic} type, whose universal cover is conformally equivalent to $S^n$,
\medskip
\myitem{(ii)} manifolds of {\bf parabolic} type, whose universal cover is conformally equivalent to $\Bbb{R}^n$, and
\medskip
\myitem{(iii)} manifolds of {\bf hyperbolic} type, consisting of all other cases.
\medskip
\noindent In the sequel, we study flat conformal structures of hyperbolic type over compact manifolds.
\medskip
\noindent Let $(\varphi,\theta)$ be a flat conformal structure over $M$. A {\bf geometric ball} in $M$ is an injective mapping 
$\alpha:B\rightarrow\tilde{M}$ from a Euclidean ball $B$ into $\tilde{M}$ such that $\varphi\circ\alpha$ is the restriction of a 
conformal mapping. Geometric balls form a partially ordered set with respect to inclusion. In \cite{Kulkarni}, it is shown that when $M$ is 
of hyperbolic type, every point of $\tilde{M}$ is contained in a maximal geometric ball. Every geometric ball carries a natural complete
hyperbolic metric. Indeed, $\partial(\varphi\circ\alpha(B))$ bounds a totally geodesic hyperplane in $\Bbb{H}^{n+1}$ and orthogonal projection defines 
a homeomorphism from $(\varphi\circ\alpha)(B)$ onto this hyperplane. The hyperbolic metric on $B$ is obtained by pulling back the metric on this
hyperplane through this orthogonal projection. We denote this metric by $g_B$. It is trivially conformal with respect to the conformal structure of $M$.
\medskip
\noindent We define the {\bf Kulkarni-Pinkall metric} $g_{KP}$ over $\tilde{M}$ by:
$$
g_{KP}(p) = \minf\left\{g_B(p)\text{ s.t. $B$ is a geometric ball and }p\in B.\right\}.
$$
\noindent This metric is the analogue in the M\"obius category of the Kobayashi metric. Trivially, $g_{KP}$ is equivariant and thus quotients to
a metric over $M$. The main result of \cite{Kulkarni} is:
\proclaim{Theorem \nextprocno\ {\bf [Kulkarni, Pinkall]}}
\noindent Let $M$ be a M\"obius manifold of hyperbolic type. Then $g_{KP}$ is positive definite and of type $C^{1,1}$.
\endproclaim
\proclabel{ThmKulkarniPinkallMetric}
\noindent Let $g_S$ be a spherical metric over $\partial_\infty\Bbb{H}^{n+1}$. Let $\overline{M}$ be the metric completion of $\tilde{M}$ with respect
to $\varphi^*g_S$. Since any two spherical metrics are uniformly equivalent, the topological space $\overline{M}$ is independant of the choice of 
spherical metric. Trivially $\varphi$ extends to a continuous map from $\overline{M}$ into $\partial_\infty\Bbb{H}^{n+1}$. We call 
$\partial\tilde{M}:=\overline{M}\setminus\tilde{M}$ the ideal boundary of $\tilde{M}$.
\medskip
\noindent Let $(B,\alpha)$ be a geometric ball. We define $C(B)$ to be the convex hull in $B$ (with respect to the hyperbolic metric) of
$\overline{\alpha(B)}\minter\partial_\infty\tilde{M}$. In proposition $4.1$ of \cite{Kulkarni}, Kulkarni and Pinkall obtain:
\proclaim{Proposition \nextprocno\ {\bf [Kulkarni, Pinkall]}}
\noindent If $M$ is a M\"obius manifold of hyperbolic type, then for every point $p\in\tilde{M}$ there exists a unique maximal
geometric ball $(B,\alpha)$ such that $p\in\alpha(C(B))$.
\endproclaim
\proclabel{PropMaximulBalls}
\noindent We denote this ball by $B(p)$. Kulkarni and Pinkall show that:
$$
g_{KP}(p) = g_{B(p)}(p).
$$
\noindent In \cite{Kulkarni}, Kulkarni and Pinkall use these maximal geometric balls to associate a canonical hyperbolic end to each flat conformal structure. These are the ends that will interest us in the sequel. We refer the reader to \cite{Kulkarni} for the details
of this construction. Let $\varphi$ be the developing map of the flat conformal structure. 
We denote the canonical hyperbolic end associated to it by $\Cal{E}(\varphi)$. Let $U\Bbb{H}^{n+1}$ be the unitary bundle of $\Bbb{H}^{n+1}$, let 
$\overrightarrow{n}:U\Bbb{H}^{n+1}\rightarrow\partial_\infty\Bbb{H}^{n+1}$ be the Gauss map and let $\pi:U\Bbb{H}^{n+1}\rightarrow\Bbb{H}^{n+1}$ be the canonical projection. 
Let $\hat{\mathi}:\tilde{M}\rightarrow U\Bbb{H}^{n+1}$ be the boundary immersion of $\Cal{E}(\varphi)$ and define $i=\pi\circ\hat{\mathi}$. $\Cal{E}(\varphi)$ has the following
useful properties:
\medskip
\myitem{(i)} $\varphi=\overrightarrow{n}\circ\hat{\mathi}$;
\medskip
\myitem{(ii)} if $p\in\tilde{M}$, if $P$ is the totally geodesic hyperplane in $\Bbb{H}^{n+1}$ normal to $\hat{\mathi}(p)$ at $i(p)$, if $g$ is the hyperbolic metric
of $P$ and if $\pi_p:\partial_\infty\Bbb{H}^{n+1}\rightarrow P$ is the orthogonal projection, then $g_{KP}(p)$ coincides with $(\pi_p\circ\varphi)^*g(p)$; and 
\medskip
\myitem{(iii)} for all $p\in\tilde{M}$, there exists a curve $\gamma:]-\epsilon,\epsilon[\rightarrow \tilde{M}$ such that $\gamma(0)=p$ and
$i\circ\gamma$ is a geodesic segment in $\Bbb{H}^{n+1}$.
\medskip
\remark Condition $(iii)$ is a strong statement about the curvature of the finite boundary of $\Cal{E}(\varphi)$, which can be defined and vanishes in the direction of the geodesic. We shall see in the sequel how this condition alone defines the geometry of the boundary immersion.
\goodbreak
\newhead{The Geodesic Boundary Property}
\noindent To better understand condition $(iii)$ of the preceeding section, we study more closely the geometry of hyperbolic ends.
\headlabel{HeadGBC}

\proclaim{Lemma \nextprocno}
\noindent Let $\tilde{N}$ be a hyperbolic end. $\tilde{N}$ is foliated by complete half-geodesics normal to the finite boundary.
\endproclaim
\proclabel{LemmaStructureFoliation}
\remark In the sequel, we will refer to this foliation as the vertical foliation.
\medskip
\proof Every subset of $\Bbb{H}^{n+1}$ is foliated in this manner. Since the transition maps preserve the concave boundary, they also preserve the foliation. The result follows.\qed
\medskip
\noindent This induces an equivalence relation on the hyperbolic end which we denote by $\sim$.
\proclaim{Lemma \nextprocno}
\noindent $\tilde{N}/\sim$ has the structure of a smooth manifold.
\endproclaim
\proclabel{LemmaStructureQuotient}
\proof Let $d$ denote the distance in $\tilde{N}$ from the finite boundary. Choose $r>0$. We claim that $N_r:=d^{-1}(\left\{r\right\})$ is a $C^{1,1}$ embedded submanifold of $\tilde{N}$. Indeed, let $\Omega\subseteq\Bbb{H}^{n+1}$ have convex boundary and let $d_\Omega$ denote the distance in $\Omega$ from the finite boundary. It follows from the properties of convex sets that $d_\Omega^{-1}(\left\{r\right\})$ is a $C^{1,1}$ embedded submanifold of $\Omega$. Since these embedded submanifolds are preserved by the transition maps, the assertion follows. Using mollifiers (c.f. \cite{SmiNLD}, for example), we obtain a smooth embedded submanifold $N'_r$ which is close to $N_r$ in the $C^1$ sense. All such embeddings have the same $C^\infty$ structure, and the result follows.\qed
\medskip
\noindent We denote $N:=\tilde{N}/\sim$.
\proclaim{Lemma \nextprocno} 
\noindent If $\tilde{N}$ is simply connected, then there exists a convex immersion, $i:N\rightarrow\Bbb{H}^{n+1}$, which is canonical up to composition by isometries of $\Bbb{H}^{n+1}$ such that:
$$
\tilde{N} = \Cal{E}(i).
$$
\endproclaim
\proclabel{LemmaStructureConstruction}
\remark In particular, if $\tilde{N}$ is an arbitrary hyperbolic end, then we may define a canonical ideal boundary $\partial_\infty\tilde{N}$ of $\tilde{N}$ as well as a canonical topology of $\tilde{N}\munion\partial_\infty\tilde{N}$.
\medskip
\proof Trivially, $N$ is simply connected. Let $d$ be the distance in $\tilde{N}$ from its finite boundary. Choose $r>0$. By the proof of Lemma \procref{LemmaStructureQuotient}, we may identify $N$ with $d^{-1}(\left\{r\right\})$. Choose $p\in N$. Let $(\alpha,U,V)$ be a coordinate chart of $\tilde{N}$ about $p$. Thus $\alpha:U\rightarrow V$, and $V \subseteq\Bbb{H}^{n+1}$ has concave boundary. Define $i_r:N\minter U\rightarrow\Bbb{H}^{n+1}$ by:
$$
i_r(q) = \alpha|_{N\minter U}.
$$
\noindent Trivially, $i_r$ is a convex immersion. Let $\hat{\mathi}_r:N\minter U\rightarrow U\Bbb{H}^{n+1}$ be the unit normal exterior vector field over $i_r$. For all $q\in N\minter U$, let $\gamma_q$ be the unit speed geodesic leaving $i_r(q)$ in the direction of $\hat{\mathi}_r(q)$. Define $\hat{\mathi}(q):N\minter U\rightarrow U\Bbb{H}^{n+1}$ by:
$$
\hat{\mathi}(q) = \partial_t\gamma_q(-r).
$$
\noindent Let $K\subseteq\Bbb{H}^{n+1}$ be a convex set such that the finite boundary of $V$ is an open subset, $\Omega$ of $\Cal{N}(K)$. Trivially, $\hat{\mathi}$ defines a homeomorphism from $N\minter U$ to $\Omega$. It follows that $\hat{\mathi}$ is a convex immersion. Moreover, $\hat{\mathi}$ is independant of $r$, and:
$$
V = \Cal{E}(\hat{\mathi}).
$$
\noindent Since $N$ is simply connected, $i_r$, $\hat{\mathi}_r$ and $\hat{\mathi}$ can be extended to mappings defined over the whole of $N$ which are canonical up to composition by homeomorphisms of $\Bbb{H}^{n+1}$. $\tilde{N}=\Cal{E}(\hat{\mathi})$, and the result follows.\qed
\medskip
\noindent The convex immersion $\hat{\mathi}:N\rightarrow\Bbb{H}^{n+1}$ yields an immersion $I:N\times]0,\infty[\rightarrow\Bbb{H}^{n+1}$ which is the end of $\hat{\mathi}$. $I$ extends continuously to a map from $N\times]0,\infty]$ to $\Bbb{H}^{n+1}\munion\partial_\infty\Bbb{H}^{n+1}$. We define $\varphi:N\rightarrow\partial_\infty\Bbb{H}^{n+1}$ by:
$$
\varphi(p) = I(p,\infty).
$$
\noindent Since $\hat{\mathi}$ is a convex immersion, $\varphi$ is a local homeomorphism. $\varphi$ thus defines a flat conformal structure over $N$. Moreover, $\varphi$ is smooth with respect to the $C^\infty$ structure of $N$. Thus the underlying $C^\infty$ structure of the flat conformal structure induced on $N$ coincides with the preexisting $C^\infty$ structure on $N$. We refer to $(N,\varphi)$ as the {\bf quotient M\"obius manifold} of the hyperbolic end $\tilde{N}$.
\medskip
\noindent Let $\tilde{N}_1$ and $\tilde{N}_2$ be hyperbolic ends. Let $(N_1,\varphi_1)$ and $(N_2,\varphi_2)$ be their respective quotient M\"obius manifolds. We define a morphism between $\tilde{N}_1$ and $\tilde{N}_2$ to be a pair $(\Phi,\tilde{\Phi})$ such that:
\medskip
\myitem{(i)} $\Phi:N_1\rightarrow N_2$ is a locally conformal mapping; 
\medskip
\myitem{(ii)} $\tilde{\Phi}:\tilde{N}_1\rightarrow\tilde{N}_2$ is a local hyperbolic isometry; and
\medskip
\myitem{(iii)} $\tilde{\Phi}$ extends to a continuous map from $\partial_\infty\tilde{N}_1=N_1$ to $\partial_\infty\tilde{N}_2=N_2$ which coincides with $\Phi$.
\medskip
\noindent In the sequel, we denote such a morphism merely by $\Phi$.
\medskip
\noindent We define a partial order ``$<$'' over the family of hyperbolic ends such that, if $\tilde{N}_1$ and $\tilde{N}_2$ are hyperbolic ends, then $\tilde{N}_1<\tilde{N}_2$ if and only if there exists an injective morphism $\tilde{\Phi}:\tilde{N}_1\rightarrow\tilde{N}_2$. If $\tilde{N}_1<\tilde{N}_2$, then we say that $\tilde{N}_1$ is contained in $\tilde{N}_2$.
\medskip
\proclaim{Definition \nextprocno, {\bf Geodesic Boundary Property}}
\noindent Let $\tilde{N}$ be a simply connected hyperbolic end. Let $N=\tilde{N}/\sim$ and let $\hat{\mathi}:N\rightarrow\Bbb{H}^{n+1}$ be the convex immersion such that $\tilde{N}=\Cal{E}(i)$. We say that $\tilde{N}$ possesses the {\bf Geodesic Boundary Property} if and only if, for every point $p\in N$ there exists:
\medskip
\myitem{(i)} a real number $\epsilon>0$;
\medskip
\myitem{(ii)} a unit speed geodesic segment $\gamma:]-\epsilon,\epsilon[\rightarrow\Bbb{H}^n$; and
\medskip
\myitem{(iii)} a continuous path $\alpha:]-\epsilon,\epsilon[\rightarrow N$,
\medskip
\noindent such that $\alpha(0)=p$ and, for all $t\in]-\epsilon,\epsilon[$:
$$
\gamma(t) = (\pi\circ\hat{\mathi}\circ\alpha)(t).
$$
\endproclaim
\proclabel{DefinitionGBC}
\remark Heuristically, $\tilde{N}$ possesses the Geodesic Boundary Property if and only if, at every boundary point, there exists a non-trivial geodesic segment passing through that point which remains in the boundary.
\medskip
\remark The Geodesic Boundary Property is a natural property of minimal convex sets in hyperbolic manifolds. Indeed, any such minimal convex set possesses the Geodesic Boundary Property, since, otherwise, there would be a point at which it would be strictly convex, and therefore be minimal.
\medskip
\remark Importantly, the Geodesic Boundary Property may be substituted by a weaker version, where, instead of a geodesic, a curve having vanishing geodesic curvature at $p$ is used. The reader may verify that this Weak Geodesic Boundary Property may be substited at every stage in the sequel where the Geodesic Boundary Property is used. As the Geodesic Boundary Property is a natural property of minimal convex sets in hyperbolic manifolds, so the Weak Geodesic Boundary Property is a natural property of minimal convex sets in more general negatively curved manifolds. We thus see how the results of this paper may be extended to a much more general setting than where they are currently presented.
\medskip
\noindent This allows us to obtain a geometric characterisation of the Kulkarni-Pinkall hyperbolic end. Let $\tilde{N}$ be a hyperbolic end. Let $d$ denote the distance in $\tilde{N}$ along the vertical foliation from the finite boundary $\partial_0\tilde{N}$ of $\tilde{N}$. For all $\delta>0$, let $N_\delta$ denote the level hypersurface $d^{-1}(\left\{\delta\right\})$. We say that $\tilde{N}$ is {\bf complete} if and only if $N_\delta$ is for some (and therefore for all) $\delta>0$.
\proclaim{Lemma \nextprocno}
\noindent Let $\tilde{N}$ be a hyperbolic end. Suppose that:
\medskip
\myitem{(i)} $\tilde{N}$ possesses the Geodesic Boundary Property; and
\medskip
\myitem{(ii)} $\tilde{N}$ is complete.
\medskip
\noindent Then $\tilde{N}$ is the Kulkarni-Pinkall hyperbolic end of its quotient flat conformal structure.
\endproclaim
\proclabel{LemmaGBCGivesKP}
\proof Let $p\in\partial_0\tilde{N}$ be a point in the finite boundary of $\tilde{N}$. Let $\msf{N}_p$ be a supporting normal to $\partial_0\tilde{N}$ at $p$ and let $H_p\subseteq\tilde{N}$ be the supporting totally geodesic hyperspace to $\partial_0\tilde{N}$ at $p$ whose normal at $p$ is $\msf{N}_p$. Since $\tilde{N}$ is complete, so is $H_p$.
\medskip
\noindent Let $K=H_p\minter\partial_0\tilde{N}$ be the intersection of $H_p$ with the finite boundary of $\tilde{N}$. Since the distance to the finite boundary in a hyperbolic end is a convex function, it follows that $K$ is a convex subset of $H_p$. Moreover, $K$ is closed. Choose $q\in K$. By the Geodesic Boundary Property, there exists $\epsilon>0$ and a unit speed geodesic segment $\gamma:]-\epsilon,\epsilon[\rightarrow H_p$ such that $\gamma(0)=q$. Since the distance to $p$ in $H_p$ is a strictly convex function, it therefore cannot attain a maximum over $K$. $K$ is therefore unbounded.
\medskip
\noindent We claim that $K$ is the convex hull of $K\minter\partial_\infty H_p$. Indeed, suppose the contrary. There exists $q\in\partial K$ which is not in the convex hull of $K\minter\partial_\infty H_p$. By rotating $\msf{N}_p$ and $H_p$ slightly around $q$, we obtain a supporting normal $\msf{N}_q$ to $\partial_0\tilde{N}$ at $q$ and a supporting totally geodesic hyperplane to $\partial_0\tilde{N}$ at $q$ whose normal is $\msf{N}_q$ such that $K\minter H_q$ is bounded. This is absurd by the previous discussion, and the assertion follows.
\medskip
\noindent It follows that $p$ is contained in the convex hull of $K\minter \partial_\infty H_p$. This condition characterises the Kulkarni-Pinkall hyperbolic end, and the result follows.\qed 
\medskip
\noindent In the compact case, moreover, the Kulkarni-Pinkall hyperbolic end is the unique maximal end. First we prove:
\proclaim{Lemma \nextprocno}
\noindent Let $\tilde{N}_1$ and $\tilde{N}_2$ be compact hyperbolic ends. Suppose, moreover that $\tilde{N}_2$ possesses the Geodesic Boundary Property. Let $(N_1,\varphi_1)$ and $(N_2,\varphi_2)$ be their respective quotient flat conformal manifolds. If $(N_1,\varphi_1)$ and $(N_2,\varphi_2)$ are isomorphic, then $\tilde{N}_1<\tilde{N}_2$. Moreover, the finite boundary, $\partial_0\tilde{N}_1$, of $\tilde{N}_1$ is a graph over the finite boundary, $\partial_0\tilde{N}_2$, of $\tilde{N}_2$.
\endproclaim
\proclabel{LemmaGBCMeansMaximality}
\proof Let $\hat{N}_1$ and $\hat{N}_2$ be the universal covers of $\tilde{N}_1$ and $\tilde{N}_2$ respectivey. Let $\hat{\Phi}_1:\hat{N}_1\rightarrow\Bbb{H}^{n+1}$ and $\hat{\Phi}_2:\hat{N}_2\rightarrow\Bbb{H}^{n+1}$ be their respective developing maps. We may assume that $\partial_\infty\hat{N}_1=\partial_\infty\hat{N}_2$ and that $\hat{\Phi}_1=\hat{\Phi}_2$ on this set.
\medskip
\noindent The identity on the ideal boundaries extends to an equivariant homeomorphism $\Psi$ from an open subset, $U_1$, of $\partial_\infty\hat{N}_1$ in $\hat{N}_1$ into an equivariant open subset, $U_2$, of $\partial_\infty\hat{N}_2$ in $\hat{N}_2$.
\medskip
\noindent Let $d:\hat{N}_1\rightarrow[0,\infty[$ be the distance in $\hat{N}_1$ to $\partial\hat{N}_1$. For all $r>0$, let $\hat{N}_{1,r}$ be the hypersurface at constant distance $r$ from $\partial\hat{N}_1$:
$$
\hat{N}_{1,r} = d^{-1}(\left\{r\right\}).
$$
\noindent For sufficiently large $r$, $\hat{N}_{1,r}$ is contained in $U$. 
\medskip
\noindent Let $V_1$ and $V_2$ be the fields of vertical vectors over $\hat{N}_1$ and $\hat{N}_2$ respectively. Let $(p_n)_\ninn\in U_1$ be a sequence converging to a point $p_0\in\partial_\infty\hat{N}_1$. Then:
$$
(\langle V_1(p_n),\Psi^*V_2(p_n)\rangle)_\ninn\rightarrow 1.
$$
\noindent Thus, by cocompactness, for sufficiently large $r$, $\Psi(\hat{N}_{1,r})$ is transverse to the field of vertical vectors over $\tilde{N}_2$. Therefore, by cocompactness, the projection from $\Psi(\hat{N}_{1,r})$ onto $\partial_0\hat{N}_2$ is a covering map, and so $\Psi(\hat{N}_{1,r})$ is a graph over $\partial_0\hat{N}_2$. Moreover, $\Psi(\hat{N}_{1,r})$ is a strict graph in the sense that it does not intersect $\partial_0\hat{N}_2$.
\medskip
\noindent Since $\hat{N}_{1,r}$ is smooth, by continuously reducing $r$, $U_1$ and $\Psi$ may be extended to contain $\hat{N}_{1,r}$ at least as long as $\Psi(\hat{N}_{1,r})$ remains a strict graph over $\partial_0\hat{N}_2$ (it will always be an immersion). Suppose therefore that there exists $r_0>0$ such that $\Psi(\hat{N}_{1,r_0})$ is not a strict graph over $\partial\hat{N}_2$ but $\Psi(\hat{N}_{1,r})$ is for all $r>r_0$.
\medskip
\noindent Suppose that $\Psi(\hat{N}_{1,r_0})$ intersects $\partial_0\hat{N}_2$ non-trivially. $\Psi(\hat{N}_{1,r_0})$ is an external tangent to $\partial\hat{N}_2$ at this point. However, by Lemma \procref{LemmaLowerCurvatureBound} the second fundamental form of $\Psi(\hat{N}_{1,r_0})$ is bounded below by $\opTanh(r_0)\opId$ in the weak sense. This therefore contradicts the Geodesic Boundary Property of $\hat{N}_2$. It follows that $\Psi(\hat{N}_{1,r_0})$ lies strictly above $\partial_0\hat{N}_2$. 
\medskip
\noindent Suppose that $\Psi(\hat{N}_{1,r_0})$ is not a graph over $\partial_0\hat{N}_2$. Then, there exists $p\in\hat{N}_{1,r_0}$ such that $\Psi(\hat{N}_{1,r_0})$ is vertical at this point. Let $q\in\partial_0\hat{N}_2$ be the vertical projection of $p$. Let $\gamma:[0,d(p,q)]\rightarrow\hat{N}_2$ be the vertical geodesic segment in $\hat{N}_2$ from $q$ to $p$. $\gamma$ lies below the graph of $\Psi(\hat{N}_{1,r})$ for all $r>r_0$. $\gamma$ is therefore an interior tangent to $\Psi(\hat{N}_{1,r_0})$ at $p$. However, as in the preceeding paragraph, $\Psi(\hat{N}_{1,r_0})$ is strictly convex at $p$, and this yields a contradiction.
\medskip
\noindent It follows that $\Psi(\hat{N}_{1,r})$ remains a strict graph over $\partial_0\hat{N}_2$ for all $r>0$. Letting $r\rightarrow 0$, it follows that $U_1=\hat{N}_{1,r}$ and that $\Psi(\partial_0\hat{N}_1)$ is a graph over $\partial_0\hat{N}_2$. The result now follows by taking quotients.\qed
\proclaim{Corollary \nextprocno}
\noindent Let $\tilde{N}$ be a compact hyperbolic end. Let $(N,\varphi)$ be its quotient flat conformal manifold. Let $\tilde{N}'$ be the Kulkarni-Pinkall hyperbolic end of $(N,\varphi)$ then $\tilde{N}$ is contained in $\tilde{N}'$ and $\partial\tilde{N}$ is a graph over $\partial\tilde{N}'$,
\endproclaim
\proclabel{CorContainedInKPEnd}
\proof The Kulkarni-Pinkall hyperbolic end satisfies the geodesic boundary condition.\qed
\proclaim{Corollary \nextprocno}
\noindent Let $\tilde{N}_1$ and $\tilde{N}_2$ be compact hyperbolic ends having the same quotient M\"obius manifold. Then there exists a unique hyperbolic end $\tilde{N}_{12}$ such that:
\medskip
\myitem{(i)} $\tilde{N}_1$ and $\tilde{N}_2$ are contained in $\tilde{N}_{12}$; and
\medskip
\myitem{(ii)} if $\tilde{N}_1$ and $\tilde{N}_2$ are contained in $\tilde{N}$, then $\tilde{N}_{12}$ is also contained in $\tilde{N}$.
\endproclaim
\proof Let $\tilde{N}_{KP}$ be the Kulkarni-Pinkall hyperbolic end of the induced flat conformal manifold. By Corollary \procref{CorContainedInKPEnd}, $\tilde{N}_1$ and $\tilde{N}_2$ are contained in $\tilde{N}_{KP}$ and $\partial_0\tilde{N}_1$ and $\partial_0\tilde{N}_2$ are graphs over $\partial_0\tilde{N}_{KP}$. Let $f_1$ and $f_2$ be their respective graph functions. The graph of $\opMin(f_1,f_2)$ in $\tilde{N}_{KP}$ is convex and yields the desired hyperbolic end.\qed
\medskip
\noindent This yields uniqueness of the maximal ends in the compact case:
\proclaim{Lemma \nextprocno}
\noindent Let $M$ be a compact M\"obius manifold. The Kulkarni-Pinkall hyperbolic end of $M$ is the unique maximal end amongst all ends whose quotient M\"obius manifold is $M$.
\endproclaim
\proclabel{LemmaKPIsUniqueMaximum}
\proof Let $\tilde{M}_{KP}$ be the Kulkarni-Pinkall hyperbolic end of M. We first show that $\tilde{M}_{KP}$ is maximal. Let $\tilde{M}$ be any other end whose quotient M\"obius manifold is $M$. Suppose that $M_{KP}<M$ and that this inclusion is strict. We thus identify $\tilde{M}_{KP}$ with a subset of $\tilde{M}$.
\medskip
\noindent Let $d$ be the distance in $\tilde{M}$ from $\partial_0\tilde{M}$. Let $p\in\partial_0\tilde{M}_{KP}$ be a point maximising distance from $\partial_0\tilde{M}$. Let $\msf{N}_p$ be a supporting normal to $\partial_0\tilde{M}_{KP}$ which is parallel to the vertical foliation of $\tilde{M}$. Let $U_p$ be the set of unit vectors, $V_p$, over $p$ in $T_p\tilde{N}$ such that:
$$
\langle V_p,\msf{N}(p)\rangle > 0.
$$ 
\noindent For all $V_p\in U$, the half geodesic in $\tilde{M}_{KP}$ leaving $p$ in the direction of $V_p$ terminates in a point in $\partial_\infty\tilde{M}_{KP}$. Let $B$ be the image of $U_p$ in $\partial_\infty\tilde{M}_{KP}$. By definition of the Kulkarni-Pinkall end, $B_p$ is a maximal ball about the image of $\msf{N}_p$.
\medskip
\noindent Let $q$ be the projection of $p\in\partial_0\tilde{M}$. Let $\msf{N}_q$ be the supporting normal to $\partial_0\tilde{M}$ at $q$ pointing towards $p$. We define $B_q$ in the same way as $B_p$. Trivially, $B_q$ contains $B_p$ in its interior, and this contradicts the maximality of $B_p$. We conclude that $\tilde{M} = \tilde{M}_{KP}$, and maximality follows.
\medskip
\noindent Let $\tilde{M}'$ be another maximal end whose quotient M\"obius manifold is $M$. Since $\tilde{M}_{KP}$ possesses the Geodesic Boundary Property, it follows by Lemma \procref{LemmaGBCMeansMaximality} that $\tilde{M}'\leqslant\tilde{M}_{KP}$. By maximality of $\tilde{M}'$, $\tilde{M}'=\tilde{M}_{KP}$, and uniqueness follows.\qed
\medskip
\noindent The proof of Theorem \procref{TheoremStructureOfKPEnd} now follows:
\medskip
{\bf\noindent Proof of Theorem \procref{TheoremStructureOfKPEnd}:\ }This follows from Lemmata \procref{LemmaGBCGivesKP} and \procref{LemmaKPIsUniqueMaximum}.\qed
\goodbreak
\newhead{The Derivative of the Curvature Operator}
\noindent Let $N$ and $M$ be Riemannian manifolds of dimensions $n$ and $(n+1)$ respectively.
The special Lagrangian curvature operator sends the space of smooth immersions from $N$ into $M$ into the space of smooth functions over $N$. These spaces may be viewed as
infinite dimensional manifolds (strictly speaking, they are the intersections of infinite nested sequences of Banach manifolds). Let $i$ be a smooth immersion from $N$ into $M$. Let ${\msf N}$ be the unit
exterior normal vector field of $i$ in $M$. We identify the space
of smooth functions over $N$ with the tangent space at $i$ of the space of smooth immersions from $N$ into $M$ as follows. Let $f:N\rightarrow\Bbb{R}$ be a smooth function. We define 
the family $(\Phi_t)_{t\in\Bbb{R}}:N\rightarrow M$ by:
\headlabel{DerivativeOfSLOp}
$$
\Phi_t(x) = \opExp(tf(x){\msf N}(x)).
$$
\noindent This defines a path in the space of smooth immersions from $N$ into $M$ such that $\Phi_0=i$. It thus defines a tangent vector to this space at $i$. Every
tangent vector to this space may be constructed in this manner.
\medskip
\noindent Let $A$ be the shape operator of $i$. This sends the space of smooth immersions from $N$ into $M$ into the space of sections of the endomorphism bundle of $TN$. We 
have the following result:
\proclaim{Lemma \nextprocno}
\noindent Suppose that $M$ is of constant sectional curvature equal to $-1$, then the derivative of the shape operator at $i$ is given by:
$$
D_i A\cdot f = f\opId -\opHess(f) - fA^2,
$$
\noindent where $\opHess(f)$ is the Hessian of $f$ with respect to the Levi-Civita covariant derivative of the metric induced over $N$ by the immersion $i$.
\endproclaim
\proclabel{DerivativeOfShapeOperator}
\proof See the proof of proposition $3.1.1$ of \cite{LabA}.\qed
\medskip
\noindent We consider the operators $\opSL_r=\opSL_r(A_\Sigma)$ and $R_\theta=R_\theta(A_\Sigma)$. Using Lemma \procref{DerivativeOfShapeOperator}, we obtain:
\proclaim{Lemma \nextprocno}
\noindent Suppose that $M$ is of constant sectional curvature equal to $-1$.
\medskip
\myitem{(i)} The derivative of $SL_r$ at $i$ is given by:
$$
(1/r)D_i SL_r\cdot f = -\opTr((\opId + r^2A^2)^{-1}\opHess(f)) + \opTr((\opId - A^2)(\opId + r^2A^2)^{-1})f.
$$
\myitem{(ii)} Likewise, the derivative of $R_\theta$ at $i$ is given by:
$$\matrix
\opTr(A(I+A^2R_\theta^2)^{-1})D_iR_\theta\cdot f \hfill&= R_\theta\opTr((\opId + r^2A^2)^{-1}\opHess(f))\hfill\cr
&\qquad\qquad + R_\theta\opTr((\opId - A^2)(\opId + r^2A^2)^{-1})f.\hfill\cr
\endmatrix$$
\endproclaim
\proclabel{DerivativeOfSLIsLaplacian}
\noindent These operators are trivially elliptic. We wish to establish when they are invertible. We first require the following technical result:
\proclaim{Lemma \nextprocno}
\noindent Let $0<n<m$ be positive integers. If $t\in ]0,\pi/2]$, then:
$$
n\opSin^2(t/n) \geqslant m\opSin^2(t/m),
$$
\noindent With equality if and only if $n=1$, $m=2$ and $t=\pi/2$.
\endproclaim
\proclabel{UsefulDecreasingFunction}
\proof The function $\opSin^2(t/2)$ is strictly convex over the interval $[0,\pi/4]$. Thus, for all $0<x<y\leqslant\pi/4$:
$$
(1/x)\opSin^2(x) < (1/y)\opSin^2(y).
$$
\noindent Thus, for $m>n\geqslant 2$, we obtain:
$$
n\opSin^2(t/n) > m\opSin^2(t/m).
$$
\noindent We treat the case $n=1$ separately. For $t\leqslant\pi/4$, the result follows as before. We therefore assume that $t>\pi/4$. Since the function 
$\opSin^2(t/2)$ is strictly concave over the interval $[\pi/4,\pi/2]$, it follows that $\opSin^2(t)\geqslant 2t/\pi$, with equality if and only 
if $t=\pi/2$. However:
$$
\opSin^2(\pi/4) = 1/2 = (2/\pi)(\pi/4).
$$
\noindent Since $m\geqslant 2$, it follows by concavity that:
$$
m\opSin^2(t/m) \leqslant \opSin^2(t),
$$
\noindent with equality if and only if $m=2$ and $t=\pi/2$. The result now follows.\qed
\medskip
\noindent Using Lagrange multipliers to determine critical points, we obtain:
\proclaim{Lemma \nextprocno}
\noindent If $\theta\geqslant (n-1)\pi/2$ and $r>\opTan(\theta/n)$, then the coefficient of the zeroth order term is non-negative:
$$
\opTr((\opId - A^2)(\opId + r^2A^2)^{-1})\geqslant 0.
$$
\noindent Moreover, this quantity reaches its minimum value of $0$ if and only if $r=\opTan(\theta/n)$ and $A$ is proportional to the identity matrix.
\endproclaim
\proclabel{ZerothOrderCoeffIsPos}
\proof For all $m$, we define the functions $\Phi_m$ and $\Theta_m$ over $\Bbb{R}^m$ by:
$$
\Phi_m(x_1,...,x_m) = \sum_{i+1}^m\frac{1-x_i^2}{1+r^2 x_i^2}, \qquad \Theta_m(x_1,...,x_m) = \sum_{i=1}^m\opArcTan(r x_i).
$$
\noindent Since the derivative of $\Theta_m$ never vanishes, $\Theta_m^{-1}(\theta)$ is a smooth submanifold of $\Bbb{R}^m$. Suppose that $\Phi_m$ achieves its minimum value on the interior of $\Theta_m^{-1}(\theta)$. Let $(\tilde{x}_1,...,\tilde{x}_m)$ be a critical point of the restriction of $\Phi_m$ to this submanifold. For all $i$, let $\tilde{\theta}_i\in [0,\pi/2[$ be such that:
$$
\opTan(\tilde{\theta}_i) = r\tilde{x}_i.
$$
\noindent Using Lagrange multipliers, we find that there exists $\eta\in[0,\pi/2]$ such that, for all $i$:
$$
\tilde{\theta}_i\in\left\{\eta,\pi/2-\eta\right\}.
$$
\noindent Let $k$ be the number of values of $i$ such that $\tilde{\theta}_i\geqslant\pi/4$. Since $\theta\geqslant(m-1)\pi/2$:
$$
k\geqslant m/2.
$$
\noindent Choose $\eta\geqslant\pi/4$. Since $\tilde{\theta}_1+...+\tilde{\theta}_m=\theta$:
$$
\eta = \frac{\theta - (m-k)\pi/2}{2k-m} = \frac{m(\theta/m) - 2(m-k)(\pi/4)}{2k-m}.
$$
\noindent If $\tilde{\Phi}_m$ is the value acheived by $\Phi_m$ at this point, then:
$$
\tilde{\Phi}_m = r^{-2}(1+r^2)(2k-m)\opCos^2(\eta) + (m-k)r^{-2}(1+r^2) - mr^{-2}.
$$
\noindent However:
$$
\pi/4 \leqslant \theta/m \leqslant \eta < \pi/2.
$$
\noindent Thus, since the function $\opCos^2$ is convex in the interval $[\pi/4,\pi/2]$:
$$
\opCos^2(\eta) \geqslant \frac{m\opCos^2(\theta/m) - 2(m-k)\opCos^2(\pi/4)}{2k-m},
$$
\noindent with equality if and only if $k=m$. Thus:
$$
\tilde{\Phi}_m \geqslant mr^{-2}(1+r^2)\opCos^2(\theta/m) - m r^{-2},
$$
\noindent with equality if and only if $\tilde{\theta}_1=...=\tilde{\theta}_m$. Since $r\geqslant\opTan(\theta/m)$, this is non-negative, and is equal to $0$ if and only if $r=\opTan(\theta/m)$. 
\medskip
\noindent We now show that $\Phi_m$ attains its minimum over $\Theta_m^{-1}(\theta)$. We treat first the case $\theta>(m-1)\pi/2$. Suppose the contrary. The functions $\Phi_m$ and $\Theta_m$ extend to continuous functions over the cube $[0,+\infty]^m$. Let $(\tilde{x}_1,...,\tilde{x}_m)$ be the point in $\Theta_m^{-1}(\theta)$ where $\Phi_m$ is minimised, and suppose now that it lies on the boundary of the cube. Since $\theta>(m-1)\pi/2$, $\tilde{x}_i>0$ for all $i$. Without loss of generality, there exists $n<m$ such that:
$$
x_1,...,x_n <+\infty,\qquad x_{n+1},...,x_m=+\infty.
$$
\noindent Let $(\tilde{\theta}_1,...,\tilde{\theta}_m)$ be as before. We define $\theta'$ by:
$$
\theta' = \tilde{\theta}_1 + ... +\tilde{\theta}_n.
$$
\noindent Since $\tilde{\theta}_{n+1}=...=\tilde{\theta}_m=\pi/2$, it follows that $\theta'=\theta - (m-n)\pi/2$. Moreover:
$$
\Phi_m(x_1,...,x_m) = \Phi_n(x_1,...,x_n) - (m-n)r^{-2}.
$$
\noindent Since $(\tilde{x}_1,...,\tilde{x}_m)$ minimises $\Phi_m$ it follows that $(\tilde{x}_1,...,\tilde{x}_n)$ is the minimal valued critical point of $\Phi_n$ in $\Theta_n^{-1}(\theta')$. Thus:
$$
\Phi_m(x_1,...,x_m) = nr^{-2}(1+r^2)\opCos^2(\theta'/n) - mr^{-2}.
$$
\noindent Let $\eta\in]0,\pi/2[$ be such that:
$$
\theta = n\pi/2 - \eta.
$$
\noindent We have:
$$
n\opCos^2(\theta'/n) = n\opSin^2(\eta/n),\qquad m\opCos^2(\theta/m) = m\opSin^2(\eta/m).
$$
\noindent It follows by Lemma \procref{UsefulDecreasingFunction} that:
$$
\Phi_m(x_1,...,x_m) > mr^{-2}(1+r^2)\opCos^2(\theta/m) - mr^{-2}.
$$
\noindent It follows that $(\tilde{x}_1,...,\tilde{x}_m)$ cannot be the minimum of $\Phi_m$ over $\Theta_m^{-1}(\theta)$, which is absurd. The result now follows in the 
case $\theta>(m-1)\pi/2$.
\medskip
\noindent It remains to study the case $\theta=(m-1)\pi/2$. This follows as before, with the single exception that it is now possible that 
$\tilde{x}_1=0$, in which case $\tilde{x}_2=...=\tilde{x}_n=+\infty$. However:
$$
\Phi_m(0,+\infty,...,+\infty) = 1 - (m-1)r^{-2}.
$$
\noindent However, $r\geqslant\opTan((m-1)\pi/2m)$. For $x\in[0,1]$, $\opTan(\pi x/4)\leqslant x$. Thus, since $m\geqslant 2$:
$$
r^{-1} \leqslant \opTan(\pi/2m) = \opTan((\pi/4)(2/m)) \leqslant 2/m.
$$
\noindent Thus:
$$
\Phi_m(0,+\infty,...,+\infty) \geqslant 1 - 4(m-1)/m^{-2} = (m-2)^2m^{-2} \geqslant 0,
$$
\noindent The result now follows.\qed
\proclaim{Lemma \nextprocno}
\myitem{(i)} If $\opSL_r(i)\geqslant(n-1)\pi/2$ and $\opTan(\opSL_r(i)/n)\leqslant r$, then $D_i\opSL_r$ is invertible.
\medskip
\myitem{(ii)} Likewise, if $\theta\geqslant(n-1)\pi/2$ and $R_\theta(i)\geqslant\opTan(\theta/n)$, then $D_iR_\theta$ is invertible.
\endproclaim
\proclabel{LemmaInvertibilityOfDerivatives}
\proof This follows immediately from the preceeding lemma, the maximum principal and the fact that second order
elliptic linear operators on the space of smooth functions over a compact manifold are Fredholm of index $0$.\qed
\goodbreak
\newhead{Deforming Equivariant Immersions}
\noindent The results of the previous section permit us to locally deform equivariant immersions of $\tilde{M}$ in $\Bbb{H}^{n+1}$. Let
$\Gamma\subseteq\opIsom(\tilde{M})$ be a cocompact subgroup acting properly discontinuously on $\tilde{M}$. Thus 
$\tilde{M}/\Gamma$ is a compact manifold. Let $\alpha:\Gamma\rightarrow\opIsom(\Bbb{H}^{n+1})$ be a homomorphism. Let 
$i:\tilde{M}\rightarrow\Bbb{H}^{n+1}$ be an immersion which is equivariant with respect to $\theta$. Thus, for all $\gamma\in\Gamma$:
\headlabel{DeformingEquivariantImmersions}
$$
i\circ\gamma = \alpha(\gamma)\circ i.
$$
\noindent Let $\rho=R_\theta(i)$. Suppose first that $i$ is an embedding. We may therefore extend $\rho$ to a smooth equivariant function over a neighbourhood of $i(\tilde{M})$ in $\Bbb{H}^{n+1}$. We obtain the following local deformation result:
\proclaim{Lemma \nextprocno}
\noindent Let $\theta\geqslant(n-1)\pi/2$ and suppose that $\rho\geqslant\opTan(\theta/n)$.
\medskip
\myitem{(i)} Let $(\alpha_t)_{t\in]-\epsilon,\epsilon[}$ be a smooth family of homomorphisms such that $\alpha_0=\alpha$;
\medskip
\myitem{(ii)} let $(\theta_t)_{t\in]-\epsilon,\epsilon[}$ be a smooth family of angles such that $\theta_0=\theta$; and
\medskip
\myitem{(iii)} let $(\rho_t)_{t\in]-\epsilon,\epsilon[}:\Bbb{H}^{n+1}\rightarrow\Bbb{R}$ be a smooth family of smooth functions such that $\rho_0=\rho$.
\medskip
\noindent There exists $0<\delta<\epsilon$ and a unique smooth family of immersions $(i_t)_{t\in]-\delta,\delta[}$ such that
$i_0=i$ and, for all $t$:
\medskip
\myitem{(i)} $R_{\theta_t}(i_t)=\rho_t\circ i_t$, and
\medskip
\myitem{(ii)} $i_t$ is equivariant with respect to $\alpha_t$.
\endproclaim
\proclabel{EquivariantDeformation}
\remark The corresponding result when $i$ is not injective is almost identical. We do not state it in order to avoid notational complexity. In the sequel, we consider embeddings inside smooth manifolds or smooth families of smooth manifolds, and so the distinction is not important.
\medskip
\proof For ease of presentation, we only prove the case where both $\rho$ and $\theta$ are constant. The general case is proven in a similar manner. The proof is divided into two stages:
\medskip
\myitem{(i)} We approximate the desired family by constructing a smooth, equivariant family of deformations of $i$ which are not necessarily immersions, 
and not necessarily of constant $\theta$-special Lagrangian curvature. First we construct a fundamental domain for $\Gamma$. Let $p$ be a point in 
$\Bbb{H}^n$. Let $P\subseteq\Bbb{H}^n$ be the orbit of $p$ under the action of $\Gamma$. Thus:
$$
P=\Gamma p.
$$
\noindent We define $\Omega\subseteq\Bbb{H}^n$ to be the set of all points on $\Bbb{H}^n$ which are closer to $p$ than to any other point in the orbit of
$p$:
$$
\Omega = \left\{q\in\Bbb{H}^n\text{ s.t. }d(q,p)<d(q,p')\text{ for all }p'\in P\setminus\left\{p\right\}\right\}.
$$
\noindent Trivially, $\Omega$ is a polyhedron and a fundemental domain for $\Gamma$.
\medskip
\noindent Using $\Omega$, we now construct the family of deformations. For each $t$, we construct a (non-continuous) deformation be defining $i_t$ to
be equal to $i$ over the interior of $\Omega$ and then extending this function to the orbit of $\Omega$ (which is almost all of $\Bbb{H}^n$) by equivariance with
respect to $\alpha_t$. These deformations may trivially be smoothed along $\partial\Omega$. The only complication is to ensure that the smoothing is performed in an equivariant
manner. The following recipe allows us to achieve exactly this.
\medskip
\noindent For any submanifold $X\in\Bbb{H}^n$ and for all $\epsilon>0$, let $X^\epsilon$ be the set of all points in $X$ which are at a distance (in $X$)
greater than $\epsilon$ from the boundary of $X$. That is:
$$
X^\epsilon = \left\{p\in X\text{ s.t. }d_X(p,\partial X)>\epsilon\right\}.
$$
\noindent Choose $\epsilon_n$ small. For all $\gamma\in\Gamma$, we define $(\tilde{\mathi}^n_t)_{t\in]\epsilon,\epsilon[}$ over
$\gamma\Omega^{\epsilon_n}$ by:
$$
\tilde{\mathi}^n_t(p) = \alpha_t(\gamma)i(\gamma^{-1}(p)).
$$
\noindent This family is trivially equivariant with respect to $(\alpha_t)_{t\in]-\epsilon,\epsilon[}$.
\medskip
\noindent Choose $\epsilon_{n-1}$ small. Let $F_{n-1}$ be any $(n-1)$-dimensional face of $\Omega$. We may trivially extend 
$(\tilde{\mathi}_t^n)_{t\in]-\epsilon,\epsilon[}$ smoothly across a neighbourhood of $F_{n-1}^{\epsilon_{n-1}}$. Since every element of $\Gamma$ is of
infinite order, there is no element which fixes any face of $\Omega$ (since otherwise it would permute the domains touching that face, and thus
be of finite order). It follows that, by choosing $\epsilon_n$ and $\epsilon_{n-1}$ small enough, we may extend this family further to a smooth
equivariant extension over every face in the orbit of $F_{n-1}$. We then continue extending this family over every face of $\Omega$ until all 
$(n-1)$-dimensional faces are exhausted. By working downwards inductively on the dimension of the faces, we thus obtain a smooth equivariant family
$(\tilde{\mathi}_t)_{t\in]-\epsilon,\epsilon[}=(\tilde{\mathi}^0_t)_{t\in]-\epsilon,\epsilon[}$ which extends $i$.
\medskip
\myitem{(ii)} We now modify this approximation to obtain the desired family of immersions. Since $\Omega$ is relatively compact, there exists 
$\delta<\epsilon$ such that, for $\left|t\right|<\delta$, $\tilde{\mathi}_t$ is an immersion. Moreover, we may suppose that for $\eta>0$ sufficiently
small, we may extend $\tilde{\mathi}_t$ smoothly along normal geodesics to a smooth equivariant immersion from $\Bbb{H}^n\times]-\eta,\eta[$ into
$\Bbb{H}^{n+1}$. We thus view $(\tilde{\mathi}_t)_{t\in]-\delta,\delta[}$ as a smooth family of immersions from $\Bbb{H}^n\times]-\eta,\eta[$ into $\Bbb{H}^{n+1}$.
\medskip
\noindent We denote by $g$ the hyperbolic metric over $\Bbb{H}^{n+1}$. We define the family $(g_t)_{t\in]-\delta,\delta[}$ such that, for all $t$:
$$
g_t = \tilde{\mathi}_t^*g.
$$
\noindent The action of $\Gamma$ over $\Bbb{H}^n$ trivially extends to an action of $\Gamma$ over $\Bbb{H}^n\times]-\eta,\eta[$. For all $t$, $g_t$ is
equivariant under this action of $\Gamma$. We denote $M=\Bbb{H}^n/\Gamma$ and we obtain a smooth family, which we also call 
$(g_t)_{t\in]-\delta,\delta[}$, of hyperbolic metrics over $M\times]-\eta,\eta[$.
\medskip
\noindent Let $j_0$ be the canonical immersion of $M$ into $M\times]-\eta,\eta[$. Trivially, with respect to $g_0$, $R_\theta(j_0)=\rho$. As in
Section \headref{DerivativeOfSLOp}, we view $R_\theta$ as a second order, non-linear differential operator sending immersions of $M$ into
$M\times]-\eta,\eta[$ into functions over $M$. Since infinitesimal variations of immersions may be interpreted as functions over $M$ times the normal
vector field of $M$ in $M\times]-\eta,\eta[$, the derivative $DR_\theta$ of $R_\theta$ may be interpreted as a second order, linear 
differential operator from $C^\infty(M)$ into $C^\infty(M)$. By Lemma \procref{LemmaInvertibilityOfDerivatives}, the operator $DR_\theta$ is invertible. After
reducing $\delta$ if necessary, the Implicit Function Theorem for non-linear PDEs therefore allows us to extend $j_0$ to a smooth family $(j_t)_{t\in]-\eta,\eta[}$ of
immersions of $M$ into $M\times]-\eta,\eta[$ such that, for all $t$, the $\theta$-special Lagrangian curvature of $j_t$ with respect to $g_t$ equals
$\rho$. For all $t$, let $\tilde{\mathj}_t$ be the lift of $j_t$ secding $\Bbb{H}^n$ into $\Bbb{H}^{n+1}$. We now define 
$i_t=\tilde{\mathi}_t\circ\tilde{\mathj}_t$. Trivially, $(i_t)_{t\in]-\delta,\delta[}$ is the desired family of immersions, and existence follows.
\medskip
\noindent Let $(i'_t)_{t\in]-\delta,\delta[}$ be another family of immersions having the desired properties. For $\delta$ sufficiently small, the image
of $i'_t$ is contained in the image of $\tilde{\mathi}_t$. For all $t$, we thus project $\tilde{\mathj}'_t = \tilde{\mathi}^{-1}_t\circ i'_t$ to an
immersion $j'_t$ of $M$ into $M\times]-\eta,\eta[$. By the uniqueness part of the Implicit Function Theorem for non-linear PDEs, for all sufficiently small
$t$, $j'_t$ coincides with $j_t$. Uniqueness now follows by a standard open/closed argument.\qed
\newhead{Compactness}
\noindent A relatively trivial variant of the reasoning used in \cite{SmiA} yields:
\proclaim{Theorem \nextprocno}
\noindent Let $M$ be a complete Riemannian manifold.
\medskip
\myitem{(i)} Let $(p_n)_\ninn,p_0\in M$ be such that $(p_n)_\ninn$ converges to $p_0$;
\medskip
\myitem{(ii)} Let $(\theta_n)_\ninn,\theta_0\in](n-1)\pi/2,n\pi/2[$ be such that $(\theta_n)_\ninn$ converges to $\theta_0$;
\medskip
\myitem{(iii)} Let $(r_n)_\ninn,r_0\in C^\infty(M)$ be strictly positive functions such that $(r_n)_\ninn$ converges to $r_0$ in the $C^\infty_\oploc$ sense; and
\medskip
\myitem{(iv)} Let $(\Sigma_n,q_n)_\ninn=(S_n,i_n,q_n)_\ninn$ be pointed, convex immersed hypersurfaces such that, for all $n$:
\medskip
\myitem{(a)} $i_n(q_n)=p_n$, and
\medskip
\myitem{(b)} $\Sigma_n$ is complete, convex and $R_{\theta_n}(i_n)=r_n\circ i_n$.
\medskip
\noindent Then, there exists a complete, pointed immersed submanifold $(\Sigma_0,q_0)=(S_0,i_0,q_0)$ in $M$ such that,
after extraction of a subsequence, $(\Sigma_n,q_n)_\ninn$ converges to $(\Sigma_0,q_0)$ in the pointed Cheeger/Gromov sense.
\endproclaim
\proclabel{CompactnessA}
\noindent The limit case where $\theta=(n-1)\pi/2$ exhibits more interesting geometric behaviour. We only require it in the constant curvature case:
\proclaim{Theorem \nextprocno}
\noindent Let $M$ be a complete Riemannian manifold.
\medskip
\myitem{(i)} Let $(p_n)_\ninn,p_0\in M$ be such that $(p_n)_\ninn$ converges to $p_0$;
\medskip
\myitem{(ii)} Let $(\theta_n)_\ninn\in[(n-1)\pi/2,n\pi/2[$ be such that $(\theta_n)_\ninn$ converges to $(n-1)\pi/2$;
\medskip
\myitem{(iii)} Let $(r_n)_\ninn,r_0\in ]0,\infty[$ be strictly positive real numbers such that $(r_n)_\ninn$ converges to $r_0$; and
\medskip
\myitem{(iv)} Let $(\Sigma_n,q_n)_\ninn=(S_n,i_n,q_n)_\ninn$ be pointed, convex immersed hypersurfaces such that, for all $n$:
\medskip
\myitem{(a)} $i_n(q_n)=p_n$, and
\medskip
\myitem{(b)} $\Sigma_n$ is convex, $R_{\theta_n}(i_n)=r_n$, and $\hat{\Sigma}_n$ is a complete submanifold of $UM$. 
\medskip
\noindent Then there exists a complete, pointed immersed submanifold $(\hat{\Sigma}_0,q_0)=(S_0,\hat{\mathi}_0,q_0)$ in 
$UM$ such that, after extraction of a subsequence, the Gauss liftings, $(\hat{\Sigma}_n,q_n)_\ninn$ converge to $(\hat{\Sigma}_0,q_0)$ in the pointed Cheeger/Gromov sense. Moreover:
\medskip
\myitem{(i)} either there exists a convex, immersed hypersurface $\Sigma_0$ in $M$ of constant $(n-1)\pi/2$-special Lagrangian curvature equal to $r_0$ such that $\hat{\Sigma}_0$ is the Gauss lifting of $\Sigma_0$ (in other words, if $\pi:UM\rightarrow M$ is the canonical projection, then $\pi\circ\hat{\mathi}_0$ is an immersion);
\medskip
\myitem{(ii)} or $\hat{\Sigma}_0$ is a covering of a complete sphere bundle over a complete geodesic.
\endproclaim
\proclabel{CompactnessB}
\remark Heuristically, if $(\Sigma_n,p_n)_\ninn=(S_n,i_n,p_n)_\ninn$ is a sequence of pointed, immersed submanifolds of constant $(n-1)\pi/2$-special Lagrangian curvature equal to $r$, then
$(\Sigma_n,p_n)_\ninn$ subconverges to $(\Sigma_0,i_0,p_0)$ where $\Sigma_0$ is either another such immersed submanifold or a complete geodesic. This (slightly abusive) language will be use in the
sequel.
\goodbreak
\newhead{The Geometric Maximum Principal}
\noindent Let $\Cal{E}$ be a hyperbolic end possessing the Geodesic Boundary Property and let $\partial_0\Cal{E}$ be its finite boundary. For all $d$, let $M_d$ be the hypersurface in $\Cal{E}$ at a distance $d$ from $\partial_0\Cal{E}$. We make the following definition:
\headlabel{UpperAndLowerBounds}
\proclaim{Definition \nextprocno}
\noindent Let $M$ be a manifold and let $\Sigma=(S,i)$ be a $C^0$ convex immersed hypersurface in $M$. Let $A$ be a family of positive definite, symmetric, bilinear forms defined on the supporting tangent planes of $\Sigma$.
The second fundamental form of $\Sigma$ at $p$ is said to be at least (resp. at most) $A$ in the weak sense if and only if, for all $p\in S$ and for each supporting tangent space $E_p$ of $\Sigma$ at $p$, there exists a smooth, convex, immersed submanifold $\Sigma'=(S,i')$ which is an exterior (resp. interior) tangent to $\Sigma$ with tangent space $E_p$ at $p$ and whose second fundamental form is bounded below (resp. above) by $A(E_p)$.
\medskip
\noindent Likewise, if $p\in S$, if $\theta\in]0,n\pi/2[$ and if $r>0$, then the $\theta$-special Lagrangian curvature of 
$\Sigma$ at $p$ is said to be at least (resp. at most) $r$ in the weak sense if and only if there exists a smooth, convex, immersed submanifold $\Sigma'=(S',i')$ of $\theta$-special Lagrangian curvature equal to $r$ which is an exterior (resp. interior) tangent to $\Sigma$ at $p$.
\endproclaim
\proclabel{DefinitionWeakCurvatureBounds}
\remark If the second fundamental form of $\Sigma$ is bounded above and below, then $\Sigma$ is necessarily of type $C^{1,1}$.
\medskip
\noindent This definition is well adapted to the Geometric Maximum Principal, whose proof requires the following result concering symmetric matrices:
\goodbreak
\proclaim{Lemma \nextprocno, {\bf Minimax Principal}}
\noindent Let $A$ be a symmetric matrix of rank $n$. If $\lambda_1\leqslant...\leqslant\lambda_n$ are the
eigenvalues of $A$ arranged in ascending order, then, for all $k$:
$$
\lambda_k = \minf_{\opDim(E)=k}\msup_{v\in E\setminus\{0\}}\langle Av,v\rangle/\|v\|^2.
$$
\endproclaim
\proof Let $e_1, ..., e_n$ be the eigenvectors of $A$. We define $\hat{E}$ by:
$$
\hat{E} = \langle e_1, ..., e_k\rangle.
$$
\noindent Let $\pi$ be the orthogonal projection onto $\hat{E}$. Let $E$ be a subspace of $\Bbb{R}^n$ of dimension $k$. For all $v$ in
$E$:
$$
\langle A\pi(v),\pi(v)\rangle\|v\|^2 \leqslant \langle Av,v\rangle\|\pi(v)\|^2.
$$
\noindent If the restriction of $\pi$ to $E$ is an isomorphism, then it follows that:
$$
\lambda_k=\msup_{v\in\hat{E}\setminus\{0\}}\langle Av,v\rangle/\|v\|^2\leqslant\msup_{v\in E\setminus\{0\}}\langle Av,v\rangle/\|v\|^2.
$$
\noindent Otherwise, there exists a non-trivial $v\in E$ such that $\pi(v)=0$, in which case:
$$
\langle Av,v \rangle \geqslant\lambda_{k+1}\|v\|^2 \geqslant\lambda_k\|v\|^2.
$$
\noindent The result now follows.\qed
\proclaim{Corollary \nextprocno}
\noindent Let $A$,$A'$ be two symmetric matrices of rank $n$ such that $A'\leqslant A$. If $\lambda_1,...,\lambda_n$ and
$\lambda_1',...,\lambda_n'$ are the eigenvalues of $A$ and $A'$ respectively arranged in ascending order, then, for all $k$:
$$
\lambda_k'\leqslant\lambda_k.
$$
\endproclaim
\noindent We now obtain the Geometric Maximum Principal for hypersurfaces of constant special Lagrangian curvature:
\proclaim{Lemma \nextprocno}
\noindent Let $M$ be a Riemannian manifold and let $\Sigma=(S,i)$ and $\Sigma'=(S',i')$ be $C^0$ convex, immersed hypersurfaces in $M$. For $\theta\in]0,n\pi/2[$, let $R_\theta$ and $R_\theta'$ be the $\theta$-special Lagrangian curvatures of $\Sigma$ and $\Sigma'$ respectively. If $p\in S$ and $p'\in S'$ are such that $q=i(p)=i'(p')$, and $\Sigma'$ is an interior tangent to $\Sigma$ at $q$, then:
$$
R_\theta(p) \geqslant R'_\theta(p').
$$
\endproclaim
\proclabel{LemmaGeomMaxPrinc}
\proof If $A$ and $A'$ are the shape operators of $\Sigma$ and $\Sigma'$ respectively, then:
$$
A'(p')\geqslant A(p).
$$
\noindent It follows that:
$$
\opArcTan(R_\theta(p)A'(p'))\geqslant\opArcTan(R_\theta(p)A(p))=\theta=\opArcTan(R_\theta'(p')A'(p')).
$$
\noindent The result now follows since the mapping $\rho\mapsto\opArcTan(\rho A'(p'))$ is strictly increasing.\qed 
\proclaim{Lemma \nextprocno}
\noindent For all $d>0$, the second fundamental form of $M_d$ is at least $\opTanh(d)\opId$ in the weak sense.
\endproclaim
\proclabel{LemmaLowerCurvatureBound}
\proof It suffices to calculate the second fundamental form of a hypersurface equidistant from a supporting totally geodesic submanifold at some point of $\partial\Cal{E}$. The result now follows from Lemma
\procref{DerivativeOfShapeOperator}.\qed
\medskip
\proclaim{Corollary \nextprocno}
\noindent Let $\theta\in]0,n\pi/2[$ be an angle. For all $d>0$, the $\theta$-special Lagrangian curvature of $M_d$ is at least $\opTan(\theta/n)/\opTanh(d)$ in the weak sense.
\endproclaim
\proclabel{CurvatureOfLevelSetsI}
\noindent For $d>0$, define the matrix $A_0(d)$ by:
$$
A_0(d) = \pmatrix \opTanh(d)\hfill& \cr &\opCoth(d)\opId_{n-1}\hfill\cr\endpmatrix,
$$
\noindent where $\opId_{n-1}$ is the $(n-1)$-dimensional identity matrix.
\proclaim{Lemma \nextprocno}
\noindent For all $d>0$, there exists a (not necessarily continuous) field $A$ of symmetric, bilinear forms over $M_d$ such that:
\medskip
\myitem{(i)} for all $p\in M_d$, $A(p)$ is conjugate to $A_0$; and
\medskip
\myitem{(ii)} the second fundamental form of $M_d$ is bounded above by $A$ in the weak sense.
\endproclaim
\proclabel{LemmaUpperCurvatureBound}
\proof For all $q\in\partial\Cal{E}$, there is a geodesic segment passing through $p$ which remains in $\partial\Cal{E}$. Thus, for all $p\in M_d$, there is a cylinder at a distance $d$ from a geodesic segment which is an interior tangent to $M_d$ at $p$. By Lemma \procref{DerivativeOfShapeOperator}, the second fundamental form of this cylinder is conjugate to $A_0$. The upper bound of the curvature at $p$ thus follows.\qed
\proclaim{Corollary \nextprocno}
\noindent Let $\theta\in](n-1)\pi/2,n\pi/2[$ be an angle. There exists a function $\kappa:[0,+\infty[\rightarrow[0,+\infty[$, which tends to $+\infty$ as $d$ tends to $0$, such that the $\theta$-special Lagrangian curvature of $M_d$ is at most 
$\kappa(d)$ in the weak sense.
\endproclaim
\proclabel{CurvatureOfLevelSetsII}
\noindent We now obtain upper and lower bounds for the distance between a hypersurface of bounded $\theta$-special Lagrangian curvature and $\partial\Cal{E}$:
\proclaim{Lemma \nextprocno}
\noindent Let $\Cal{E}$ be a hyperbolic end. Let $\partial\Cal{E}$ be the boundary of $\Cal{E}$. 
Let $\theta\in](n-1)\pi/2,n\pi/2[$ be an angle. There exists a decreasing function $\delta:[\opTan(\theta/n),+\infty[\rightarrow]0,+\infty[$ such that if $r\leqslant R\in]\opTan(\theta/n),\infty[$ and if $\Sigma=(S,i)$ is a compact, 
convex immersed submanifold such that $R_\theta(i)\in[r,R]$, then, for all $p\in S$:
$$
\delta(R) \leqslant d(i(p),\partial\Cal{E}) \leqslant \opArcTanh(r^{-1}\opTan(\theta/n)).
$$
\endproclaim
\proclabel{UpperBoundInHyperbolicEnd}
\proof For all $\rho>0$, let $M_\rho$ be the level hypersurface in $\Cal{E}$ at a distance of $R$ from $\partial\Cal{E}$. Since $\Sigma$ is compact, there exists a point $p\in S$ maximising the distance from $\partial\Cal{E}$. Let $d$ be the distance of $i(p)$ from $\partial\Cal{E}$. $\Sigma$ is an interior tangent to $M_d$ at $p$. The upper bound now follows by Lemma \procref{CurvatureOfLevelSetsI} and the geometric maximum principle (Lemma \procref{LemmaGeomMaxPrinc}). The 
lower bound follows in an analogous way, using Lemma \procref{CurvatureOfLevelSetsII} instead of Lemma \procref{CurvatureOfLevelSetsI}.\qed
\goodbreak
\newhead{Uniqueness}
\noindent We show that the metric induced by $i$ is uniformly equivalent, up to reparametrisation, with the Kulkarni-Pinkall metric:
\goodbreak
\proclaim{Lemma \nextprocno}
\noindent Let $\theta\in](n-1)\pi/2,n\pi/2[$ be an angle, and let $r>\opTan(\theta/n)$ be a positive real number. There exists $K = K(r,\theta, n)>0$
which only depends on $r$, $\theta$ and $n$ such that:
\medskip
\myitem{(i)} if $M$ is a compact manifold and $(\varphi,\theta)$ is a flat conformal structure of hyperbolic type over $M$;
\medskip
\myitem{(ii)} if $i:M\rightarrow\Bbb{H}^{n+1}$ is a complete, equivariant, convex immersion such that $R_\theta(i)=r$ and $\overrightarrow{n}\circ\hat{\mathi}=\varphi$; and
\medskip
\myitem{(iii)} if $\alpha:M\rightarrow M$ is a reparametrisation such that $i\circ\alpha$ is a graph over $\hat{\mathj}$, where $\hat{\mathj}$ is the boundary
immersion of $\Cal{E}(\varphi)$,
\medskip
\noindent then, if $g$ is the hyperbolic metric on $\Bbb{H}^{n+1}$:
$$
K^{-1} g_{KP} \leqslant (i\circ\alpha)^*g\leqslant Kg_{KP}.
$$
\endproclaim
\proclabel{LemmaUnifEquiv}
\proof Let $\Cal{E}(\varphi)$ be the Kulkarni-Pinkall hyperbolic end of $\varphi$. Since, in particular, $i$ is a convex immersion, by Lemma \procref{LemmaGBCMeansMaximality}, $i$ may be viewed as an immersion
from $M$ into $\Cal{E}(\varphi)$. For all $R>0$, let $M_R$ be the hypersurface at distance $R$ from $\partial_0\Cal{E}(\varphi)$. By Lemma 
\procref{UpperBoundInHyperbolicEnd}, there exists $R>\epsilon>0$ such that $i(M)$ lies between $M_\epsilon$ and $M_R$. Define $\pi:M\rightarrow\partial_0\Cal{E}(\varphi)$
such that $\pi(p)$ is the orthogonal projection of $i(p)$ onto $\partial_0\Cal{E}(\varphi)$. For all $p\in M$, let $\gamma_p$ be the geodesic segment
joining $\pi(p)$ to $i(p)$. Let $\msf{N}_p$ be the exterior normal to $i(M)$ at $p$.
\medskip
\noindent We show that there exists $\delta$, which only depends on $r$, $\theta$ and $n$ such that $\gamma_p$ makes an angle of at most $\pi/2-\delta$
with $\msf{N}_p$. We consider the universal covers of $M$ and $\Cal{E}(\varphi)$. In this case $i(M)$ only intersects $\gamma_p$ once in $B_\epsilon(i(p))$. Let $(M_n,p_n)_\ninn$ be
a sequence of complete, pointed manifolds. For all $n$, let $(\theta_n,\varphi_n)$ be a flat conformal structure of hyperbolic type over $M_n$ and let $i_n:\tilde{M}_n\rightarrow\Bbb{H}^{n+1}$ be 
a complete, equivariant, convex immersion such that $\varphi_n=\overrightarrow{n}\circ\hat{\mathi}_n$. For all $n$, let $\gamma_n$ be the geodesic segment joining $\pi_n(p_n)$ to $i_n(p_n)$.
Suppose that the angle that $\gamma_n$ makes with $\msf{N}_{p_n}$ tends to $\pi/2$. 
\medskip
\noindent By Theorems \procref{CompactnessA} and \procref{CompactnessB}, after extracting a subsequence, we may assume that the sequences $(i_n,M_n,p_n)_\ninn$ and $(\gamma_n)_\ninn$ converge to $(i_0,M_0,p_0)$ and $\gamma_0$
respectively. The limit hypersurface $(i_0,M_0,p_0)$ is an immersed submanifold in $\Bbb{H}^{n+1}$. Since the
$\gamma_n$ have length bounded below by $\epsilon$, $\gamma_0$ is a finite length geodesic segment which is an interior tangent to $i_0$ at $p_0$. This is impossible,
since $i_0$ is strictly convex, and the result follows.
\medskip
\noindent For $p\in M$, let $P_p$ be the supporting totally geodesic hyperspace to $\Cal{E}(\varphi)$ normal to $\gamma_p$ at $\pi(p)$. Since $i(M)$ lies below $M_R$ and since its normal makes an angle of at most $\pi/2-\delta$ with $\gamma_p$, there exists $K$, which only depends on $R$,$\epsilon$ and $\delta$ such that the normal projection from $i(M)$ onto $P_p$ is $K$-bilipschitz at $p$. The result now follows by the relationship between
$\Cal{E}(\varphi)$ and $g_{KP}$.\qed
\medskip
\noindent This yields uniqueness:
\proclaim{Lemma \nextprocno\ {\bf Uniqueness}}
\noindent Let $M$ be a conformally flat manifold of hyperbolic type. Let $\alpha:\pi_1(M)\rightarrow\opIsom(\Bbb{H}^{n+1})$ be the holonomy and let $\varphi:\tilde{M}\rightarrow\partial_\infty\Bbb{H}^{n+1}$ be the developing map.
\medskip
\noindent Let $\theta\in[(n-1)\pi/2,n\pi/2[$ be an angle, and let $r\geqslant\opTan(\theta/n)$. Let $i,i':\tilde{M}\rightarrow\Bbb{H}^{n+1}$ be complete, equivariant, convex immersions such that $R_\theta(i)=R_\theta(i')=r$ and $\overrightarrow{n}\circ\hat{\mathi}=\overrightarrow{n}\circ\hat{\mathi}'=\varphi$. Then, up to reparametrisation, $i=i'$.
\medskip
\noindent Moreover $i=i'$ is a graph over the finite boundary of the Kulkarni-Pinkall hyperbolic end of $M$, and is thus strictly contained within this hyperbolic end.
\endproclaim
\proclabel{LemmaUniqueness}
\proof By Lemma \procref{LemmaGBCMeansMaximality}, we view $i$ and $i'$ as immersions inside $\Cal{E}(\varphi)$. We first consider the case where
$\theta\neq(n-1)\pi/2$ and extend $i$ and $i'$ to unique foliations $(i_t)_{t\in[r,+\infty[}$ and $(i'_t)_{t\in[r,+\infty[}$ respectively which cover the lower end 
of $\Cal{E}(\varphi)$.
\medskip
\noindent Let $I\subseteq [r,+\infty[$ be such that, for all $T\in I$, there exists a foliation $(i^T_t)_{t\in [r,T[}$ of $\Cal{E}(\varphi)$ such that
$i_r=i$ and, for all $t$, $R_\theta(i_t)=t$. By the local uniqueness part of Lemma \procref{EquivariantDeformation}, these foliations are unique. In
other words, for all $r\leqslant t<T<T'$:
$$
i_t^T = t_t^{T'}.
$$
\noindent By Lemma \procref{EquivariantDeformation}, there exists $\delta>0$ and a smooth family $(i_t)_{t\in[r,r+\delta[}$ such that $i_r=r$, and, for all $t$, 
$R_\theta(i_t)=t$. Let $\msf{N}$ be the normal vector field over $i$. Let $f$ be the function over $M$ such that $f\msf{N}$ is the infinitesimal deformation
of $(i_t)_{t\in[r,r+\delta[}$. Then:
$$
D_iR_\theta f = 1\geqslant 0.
$$
\noindent It follows by Lemma \procref{ZerothOrderCoeffIsPos} that $f<0$. Thus,
by reducing $\delta$ if necessary, $(i_t)_{t\in[r,r+\delta[}$ is a foliation. $I$ is therefore non-empty. Let $T$ be the suprememum of $I$ and suppose that
$T<+\infty$. By uniqueness, there exists a foliation $(i_t)_{t\in [r,T[}$ with the given properties.
\medskip
\noindent For all $t\in [r,T[$, by Lemma \procref{LemmaGBCMeansMaximality}, $i_t$ is a graph over $\partial\Cal{E}(\varphi)$. Since
$(i_t)_{t\in [r,T[}$ is a foliation, the corresponding graphs form a monotone family. In fact, the graphs are monotone decreasing. For all $t$, let 
$\opVol_t$ and $\opInj_t$ be the volume and injectivity radius respectively of $i_t$. By Lemma \procref{LemmaUnifEquiv}, $\opVol_t$ is uniformly bounded above and
$\opInj_t$ is uniformly bounded below for $t\in T$. It follows by Theorem \procref{CompactnessA} that, for every sequence 
$(t_n)_\ninn$ which converges to $T$, $(i_{t_n})_\ninn$ subconverges. By monotonicity, all these subsequences converge to the same immersion, and thus $(i_t)_{t\in[r,T]}$ converges
as $t$ tends to $T$. We thus extend $(i_t)_{t\in[r,T[}$ to a foliation $(i_t)_{t\in [r,T]}$ defined over the closed interval.
\medskip
\noindent Applying Lemma \procref{EquivariantDeformation} again, this foliation can be extended to a foliation $(i_T)_{t\in[r,T+\delta[}$. This contradicts the definition of
$T$. We thus obtain the desired foliation.
\medskip
\noindent Let $f$ and $f'$ be the functions of which $i$ and $i'$ are the graphs over $\partial\Cal{E}(\varphi)$. Suppose that $f'<f$ at some point. For all $R$, let 
$M_R$ be the hypersurface of $\Cal{E}(\varphi)$ at distance $R$ from $\partial\Cal{E}(\varphi)$. Let $\epsilon>0$ be such that $i$ and $i'$ lie above $M_\epsilon$. By 
Lemma \procref{UpperBoundInHyperbolicEnd}, $(i_t)_{t\in[r,+\infty[}$ converges to $\partial\Cal{E}(\varphi)$ in the Hausdorff sense as $t$ tends to $+\infty$. In particular, there exists $R_0>r$ 
such that $i_R$ lies below $M_\epsilon$ and thus does not intersect $i'$. Let $R$ be the supremum of all $s\in[r,R_0]$ such that $i_s$ intersects $i'$ non-trivially. By 
compactness $i_R$ is an interior tangent to $i'$ at some point. However, $R_\theta(i_r)=R>R_\theta(i')$, which is a contradiction by the Geometric Maximum Principal
(Lemma \procref{LemmaGeomMaxPrinc}).
\medskip
\noindent It follows that $f'\geqslant f$. By symmetry, $f\geqslant f'$, and the result now follows for $\theta\neq(n-1)\pi/2$.
\medskip
\noindent Suppose that $\theta=(n-1)\pi/2$. By Lemma \procref{EquivariantDeformation}, there exist smooth families $(i_\eta)$ and $(i'_\eta)$ for $\eta\in[(n-1)\pi/2,(n-1)\pi/2+\delta[$ 
such that $i=i_{(n-1)\pi/2}$, $i'=i'_{(n-1)\pi/2}$ and, for all $\eta$:
$$
R_\eta(i_\eta) = R_\eta(i'_\eta) = r.
$$
\noindent By uniqueness for $\theta\neq(n-1)\pi/2$, $i_\eta=i'_\eta$ for all $\eta\neq(n-1)\pi/2$ and the result now follows for $\theta=(n-1)\pi/2$ by taking
limits.\qed
\goodbreak
\newhead{Main Results}
\noindent Existence follows from Theorem $1.4$ of \cite{SmiNLD}. For the reader's convenience, we include a proof based on the more elementary Theorem $1.2$ of the same paper. Throughout the rest of this section, a convex set will be said to be $\epsilon$-convex for some $\epsilon>0$ if and only if its second fundamental form with respect to every supporting normal is bounded below by $\epsilon\opId$ in the weak sense. We quote Theorem $1.2$ of \cite{SmiNLD}:
\proclaim{Theorem \nextprocno}
\noindent Choose $\theta\in[(n-1)\pi/2,n\pi/2[$. Let $H\subseteq\Bbb{H}^{n+1}$ be a totally geodesic hypersurface. Let $\Omega\subseteq H$ be a bounded open subset. Let $\hat{\Sigma}\subseteq\Bbb{H}^{n+1}$ be a convex hypersurface which is a graph over $\Omega$ such that $\partial\hat{\Sigma} = \partial\Omega$ and:
$$
R_\theta(\hat{\Sigma}) \leqslant R_1,
$$
\noindent in the weak sense, where $R_1\geqslant \opTan^{-1}(\theta/n)$. If $\theta>(n-1)\pi/2$, then, for all $r\in[R_1,\infty]$, there exists a unique immersed hypersurface $\Sigma_r\subseteq\Bbb{H}^{n+1}$ such that:
\medskip
\myitem{(i)} $\Sigma_r$ is $C^0$ and $C^\infty$ in its interior;
\medskip
\myitem{(ii)} $\partial\Sigma_r = \partial\Omega$;
\medskip
\myitem{(iii)} $\Sigma_r$ is a graph over $\Omega$ lying below $\hat{\Sigma}$; and
\medskip
\myitem{(iv)} $R_\theta(\Sigma_r) = r$.
\medskip
\noindent Moreover, the same result holds for $\theta=(n-1)\pi/2$ provided that, in addition, $\hat{\Sigma}$ is $\epsilon$-convex, for some $\epsilon>0$.
\endproclaim
\proclabel{TheoremLocalExistence}
\remark The statement of this theorem differs slightly from that appearing in \cite{SmiNLD} because (for technical reasons) the special Lagrangian curvature as defined in \cite{SmiNLD} is the reciprocal of the special Lagrangian curvature as defined here.
\medskip
\noindent Following \cite{GuanSpruck} and \cite{SmiNLD}, we use the Perron method to obtain:
\proclaim{Lemma \nextprocno}
\noindent Let $\Cal{E}$ be a hyperbolic end satisfying the Geodesic Boundary Condition. For all $\theta\in](n-1)\pi/2,n\pi/2[$ and for all $r>\opTan(\theta/n)$, there exists a strictly convex immersed hypersurface $\Sigma=(S,i)$ in $\Cal{E}$ which is a graph over the finite boundary of $\Cal{E}$ such that $R_\theta(i)=r$.
\medskip
\noindent Moreover, if the quotient M\"obius manifold is neither $S^{n-1}\times S^1$ nor $S^{n-1}\times\Bbb{R}$, where $S^k$ is the $k$-dimensional sphere, then the same result holds for $\theta=(n-1)\pi/2$.
\endproclaim
\proclabel{LemmaExistence}
\proof We first treat the case where the quotient M\"obius manifold of $\Cal{E}$ is compact and $\theta>(n-1)\pi/2$. Let $\partial_0\Cal{E}$ be the finite boundary of $\Cal{E}$. For $d>0$, let $\Sigma^0_d$ be the level hypersurface at distance $d$ from $\partial_0\Cal{E}$. By Lemma \procref{LemmaLowerCurvatureBound}, the second fundamental form of $\Sigma^0_d$ is greater than $\opTanh(d)\opId$ in the weak sense. Since $\opTanh(d)$ tends to $1$ as $d$ tends to $+\infty$, for sufficiently large $d$, the $\theta$-special Lagrangian curvature of $\Sigma^0_d$ is at most $r$ in the weak sense. Choose such a $d$ and denote $\Sigma_0 = \Sigma^0_d$.
\medskip
\noindent By definition, $\Sigma_0$ is a graph over $\partial_0\Cal{E}$. Let $f_0$ be the function whose graph $\Sigma_0$ is. Let $\Sigma_1$ be a strict graph over $\partial_0\Cal{E}$ lying below $\Sigma_0$ such that $R_\theta(\Sigma_1)\leqslant r$ in the weak sense. There exists $\epsilon>0$, which only depends on $\theta$ and $r$ such that $\Sigma_1$ is $\epsilon$-convex. In particular, by Lemma \procref{LemmaUpperCurvatureBound} and the Geometric Maximum Principal, there exists $\delta>0$ such that $\Sigma_1$ lies at a distance of at least $\delta$ from $\partial_0\Cal{E}$. Let $U_1$ be the open set lying between $\partial_0\Cal{E}$ and $\Sigma_1$. Choose $p\in\Sigma_1$. Let $\msf{N}_p$ be a supporting normal to $\Sigma_1$ at $p$ chosen such that, for any other supporting normal $\msf{N}_p'$ to $\Sigma_1$ at $p$:
$$
\langle\msf{N}'_p,\msf{N}_p\rangle \geqslant \eta,
$$
\noindent for some $\eta>0$. Such an $\msf{N}_p$ always exists since $\Sigma_1$ bounds a strictly convex set with non-trivial interior (c.f. Lemma $4.7$ of \cite{SmiNLD}). Let $\delta_1>0$ be smaller than the injectivity radius of $\Cal{E}$ at $p$. Let $\gamma$ be the unit speed geodesic such that:
$$
\partial_t\gamma(0) = \msf{N}_p.
$$
\noindent For small $t$, let $D_{p,t}$ be the totally geodesic disk in $\Cal{E}$ of radius $\delta_1$ about $\gamma(t)$ whose exterior normal at $p$ is $\partial_t\gamma(t)$. By strict convexity, $D_{p,0}$ only intersects $\Sigma_1$ at a single point. There therefore exists $\delta_2>0$ such that, for all $t\in ]-\delta_2,0[$, $\Omega_t:=U_0\minter D_{p,t}$ is a convex set and the portion of $\Sigma_1$ lying above $\Omega_t$ is a graph over $\Omega_t$ which we denote by $\Sigma_{1,t}$. Moreover, $\delta_2$ may also be chosen sufficiently small such that it doesn't intersect $\partial_0\Cal{E}$.
\medskip
\noindent By Theorem \procref{TheoremLocalExistence}, for all $t\in]-\delta_2,0[$, there exists a unique graph $\Sigma'_{1,t}$ over $\Omega_t$, lying beneath $\Sigma_{1,t}$ such that:
$$
R_\theta(\Sigma_{1,t}') = r.
$$
\noindent For all $t\in]-\delta_2,0[$, let $\Sigma'_t$ be the hypersurface obtained by replacing the portion $\Sigma_{1,t}$ of $\Sigma_1$ with $\Sigma'_{1,t}$. By uniqueness, this is a continuous family. Moreover, for $t_1>t_2$, $\Sigma_{t_1}'$ lies above $\Sigma_{t_2}'$.
\medskip
\noindent We claim that $R_\theta(\Sigma'_t)\leqslant r$ in the weak sense. It suffices to verify this property along $\partial\Omega=\partial\Sigma'_{1,t}$. However, along $\partial\Omega$, this property follows by the convexity of the curvature condition ($R_\theta$ is a convex function, c.f. Lemma $2.4$ of \cite{SmiNLD}). The assertion therefore follows.
\medskip
\noindent In particular, $\Sigma'_t$ is $\epsilon$-convex for all $t$. We claim that $\Sigma'_t$ is a graph over $\partial_0\Cal{E}$. Indeed, since $D_{p,t}$ lies strictly above $\partial_0\Cal{E}$, so does $\Sigma_t'$ for all $t$. $\Sigma_t'$ therefore only ceases to be a graph if it becomes vertical at some point $q_0$ for some value $t_0$ of $t$. $t_0$ may be chosen such that $\Sigma_t$ is a graph over $\partial_0\Cal{E}$ for all $t\in]t_0,0[$. Let $\underline{q}_0$ be the projection of $q$ in $\partial_0\Cal{E}$. Let $\gamma:[0,d(\underline{q}_0,q_0)]\rightarrow\Cal{E}$ be the geodesic segment in $\Cal{E}$ joining $\underline{q}_0$ to $q_0$. For all $t$, let $U'_t$ be the open set lying between $\partial_0\Cal{E}$ and $\Sigma_t'$. For $t>t_0$, since $\Sigma_t'$ lies above $\Sigma_{t_0}'$, $\gamma$ is contained in $U_t'$. It follows by continuity that $\gamma$ is contained in $U'_{t_0}$, and thus $\partial_t\gamma$ is an interior tangent to $\Sigma_{t_0}$ at $q_0$, which contradicts strict convexity. The assertion follows. 
\medskip
\noindent We choose any $t\in[-\delta_2,0]$ and define $\Sigma_2=\Sigma_t'$. We denote by $A$ this operation for obtaining new immersed hypersurfaces out of old ones. Let $\Sigma_1$ and $\Sigma_2$ be two graphs over $\partial_0\Cal{E}$ and let $f_1$ and $f_2$ be the respective functions whose graphs they are. Suppose that:
\medskip
\myitem{(i)} $f_1,f_2\leqslant f_0$; and
\medskip
\myitem{(ii)} $R_\theta(\Sigma_1),R_\theta(\Sigma_2)\leqslant r$ in the weak sense.
\medskip
\noindent Define $f_{1,2}$ by:
$$
f_{1,2} = \opMin(f_1,f_2).
$$
\noindent Let $\Sigma_{1,2}$ be the graph of $f_{1,2}$. Then $\Sigma_{1,2}$ lies below $\Sigma_0$, and, by convexity of the curvature condition (c.f. Lemma $2.4$ of \cite{SmiNLD}):
$$
R_\theta(\Sigma_{1,2}) \leqslant r.
$$
\noindent We denote this operation by $B$.
\medskip
\noindent Let $\Cal{F}$ be the family of immersed hypersurfaces in $\Cal{E}$ obtained from $\Sigma_0$ by a finite number of combinations of the operations $A$ and $B$. For any $\Sigma\in\Cal{F}$, let $f(\Sigma)$ be the function of which $\Sigma$ is the graph, and let $U(\Sigma)$ be the open set contained between $\partial_0\Cal{E}$ and $\Sigma$. Define $V_0\geqslant 0$ by:
$$
V_0 = \minf\left\{\opVol(U(\Sigma))\text{ s.t. }\Sigma\in\Cal{F}\right\}.
$$
\noindent There exists a sequence $(\Sigma_n)_\ninn\in\Cal{F}$ such that:
\medskip
\myitem{(i)} for all $n\geqslant m$:
$$
f(\Sigma_n)\leqslant f(\Sigma_m)\text{; and}
$$
\myitem{(ii)} $(\opVol(U(\Sigma_n)))_\ninn$ tends to $V_0$.
\medskip
\noindent Let $f_\infty$ be the function to which $(f(\Sigma_n))_\ninn$ converges pointwise. By Lemma \procref{LemmaUpperCurvatureBound} and the Geometric Maximum Principal, there exists $d_0>0$ such that, for all $n$:
$$
f(\Sigma_n)\geqslant d_0.
$$
\noindent It follows that $f_\infty\geqslant d_0$. Moreover, since the graphs $(f(\Sigma_n))_\ninn$ form the boundaries of a nested sequence of $\epsilon$-convex sets, the graph of $f_\infty$ is also the boundary of an $\epsilon$-convex set, and, by strict convexity as before, the graph of $f_\infty$ is never vertical. It follows that $f_\infty$ is $C^{0,1}$ and that $(f(\Sigma_n))_\ninn$ converges to $f_\infty$ in the $C^{0,\alpha}$ sense for all $\alpha$.
\medskip
\noindent We claim that $f_\infty$ is smooth. Let $\Sigma_\infty$ be the graph of $f_\infty$. Choose $p\in\Sigma_\infty$. Let $\msf{N}_p$ be a supporting normal to $\Sigma_\infty$ at $p$ chosen such that, for any other supporting normal $\msf{N}_p'$ to $\Sigma_\infty$ at $p$:
$$
\langle\msf{N}'_p,\msf{N}_p\rangle \geqslant \eta,
$$
\noindent for some $\eta>0$. Let $\delta_1>0$ be smaller than the injectivity radius of $\Cal{E}$ at $p$. Let $\gamma$ be the unit speed geodesic such that:
$$
\partial_t\gamma(0) = \msf{N}_p.
$$
\noindent For small $t$, let $D_{p,t}$ be the totally geodesic disk in $\Cal{E}$ of radius $\delta_1$ about $\gamma(t)$ whose exterior normal at $p$ is $\partial_t\gamma(t)$. By strict convexity, $D_{p,0}$ only intersects $\Sigma_\infty$ at a single point. There therefore exists $\delta_2>0$ such that, for all $t\in ]-\delta_2,0[$, $\Omega_t:=U(\Sigma_\infty)\minter D_{p,t}$ is a convex set and the portion of $\Sigma_\infty$ lying above $\Omega_t$ is a graph over $\Omega_t$.
\noindent By reducing $\delta_2$ if necessary, there exists $N\in\Bbb{N}$ such that, for all $n\geqslant N$, and for all $t\in ]-\delta_2,0[$, $\Omega_{n,t}:=U(\Sigma_n)\minter D_{p,t}$ is a convex set and the portion of $\Sigma_n$ lying above $\Omega_{n,t}$ is a graph over $\Omega_{n,t}$. Choose $t\in]-\delta_2,0[$ and for all $n\geqslant N$, define $\Sigma_n'$ by replacing the portion of $\Sigma_n$ lying above $\Omega_{n,t}$ with the smooth graph obtained from Theorem \procref{TheoremLocalExistence}. 
\medskip
\noindent $(f(\Sigma_n'))_\ninn$ is a decreasing sequence and therefore tends towards a $C^{0,1}$ limit, $f'_\infty$ in the $C^{0,\alpha}$ sense for all $\alpha$. For all $n\geqslant N$, $\Sigma_n'$ lies below $\Sigma_n$. Therefore:
$$
f'_\infty \leqslant f_\infty.
$$
\noindent We claim that $f'_\infty=f_\infty$. Indeed, suppose that $f'_\infty < f_\infty$, then:
$$
\opVol(U(f'_\infty)) < \opVol(U(f_\infty)),
$$
\noindent which contradicts the minimality of the volume below $f_\infty$. By Theorem  \procref{CompactnessA}, the portion of $(\Sigma'_n)_\ninn$ lying above $\Omega_{n,t}$ converges in the $C^\infty_\oploc$ sense to the portion of $\Sigma_\infty$ lying above $\Omega_{\infty,t}$, which is a non-trivial neighbourhood of $p$. It follows that $\Sigma_\infty$ is smooth at $p$ and that $R_\theta(\Sigma_\infty)=r$ near $p$. Since $p\in\Sigma_\infty$ is arbitrary, the result follows.
\medskip
\noindent Suppose that $\theta=(n-1)\pi/2$. Let $(\theta_n)_\ninn\in](n-1)\pi/2,n\pi/2[$ be a decreasing sequence converging towards $\theta$. Suppose moreover, that for all $n$:
$$
r> \opTan^{-1}(\theta_n/n).
$$
\noindent For all $n$, let $\Sigma_n$ be the immersed hypersurface such that:
$$
R_{\theta_n}(\Sigma_n) = r.
$$
\noindent For all $n$, let $f_n$ be the function of which $\Sigma_n$ is the graph and let $U_n$ be the open convex set lying between $\partial_0\Cal{E}$ and $\Sigma_n$. For all $d>0$, let $M_d$ be the level hypersurface at distance $d$ from $\partial_0\Cal{E}$. By Lemma \procref{LemmaLowerCurvatureBound}, there exists $D>0$ such that, for all $n$, and for all $d\geqslant D$, $R_{\theta_n}(M_d)$ is not greater than $r$. It follows by the Geometric Maximum Principal that, for all $n$, $\Sigma_n$ lies below $M_D$. There therefore exists a convex set $U_\infty$, lying below $M_D$ to which $(U_n)_\ninn$ subconverges in the Haussdorf sense.
\medskip
\noindent Let $V$ be the unit tangent vector field to the vertical foliation of $\Cal{E}$. For all $n$, since $\Sigma_n$ is a graph over $\partial_0\Cal{E}$, if $\msf{N}_n$ is the outward unit normal vector to $\Sigma_n$, then:
$$
\langle V,\msf{N}_n\rangle > 0.
$$
\noindent Taking limits, if $\msf{N}_\infty$ is a supporting normal to $U_\infty$, then:
$$
\langle V,\msf{N}_n\rangle \geqslant 0.
$$
\noindent By Theorem \procref{CompactnessB}, the sequence $(\Sigma_n)$ can only degenerate by converging towards a complete geodesic. If this happens, then the above condition on the supporting normal to $U_0$ implies one of two possibilities:
\medskip
\myitem{(i)} either this geodesic is vertical, which is impossible, since $\Sigma_n$ lies below $M_D$ for all $n$;
\medskip
\myitem{(ii)} or this geodesic coincides with $\partial_0\Cal{E}$, which is excluded by the hypotheses on $\Cal{E}$.
\medskip
\noindent We thus conclude that $\Sigma_n$ never degenerates. It follows that the boundary of $U_\infty$ is smooth. Moreover, as before, it is always transverse to $V$. It follows that $(f_n)_\ninn$ is equicontinuous, and therefore subconverges to a function, $f_\infty$. Since the graph of $f_\infty$ is the boundary of $U_\infty$, $f_\infty$ is smooth and its graph has constant $\theta$-special Lagrangian curvature equal to $r$. The concludes the proof when the quotient M\"obius manifold is compact.
\medskip
\noindent To conclude, we outline the proof in the case when the quotient M\"obius manifold is not compact. Let $(U_n)_\ninn$ be an exhaustion of $\partial_0\Cal{E}$ by relatively compact open sets. For each $n$, we verify that the Perron method preserves graphs over $U_n$, and we thus obtain a smooth graph over $U_n$ of constant special Lagrangian curvature. Moreover, using the Geometric Maximum Principal, we show that these graphs are uniformly bounded, and thus subconverge to a smooth graph over the whole of $\partial_0\Cal{E}$ which has the desired properties. The general result now follows.\qed
\medskip
{\bf\noindent Proof of Theorem \procref{TheoremExistenceAndUniqueness}:\ }This is the union of Lemmata \procref{LemmaUniqueness} and \procref{LemmaExistence}.\qed
\medskip
{\bf\noindent Proof of Theorem \procref{TheoremFoliations}:\ }Using Lemma \procref{EquivariantDeformation}, these hypersurfaces form a smooth family. Moreover, we can show that the derivative of
$i_{r,\theta}$ with respect to $r$ is strictly negative. Thus, if $r'<r$ are close, then $\Sigma_{r,\theta}$ lies 
strictly below $\Sigma_{r',\theta}$. It follows that this family defines a foliation. By Lemma 
\procref{UpperBoundInHyperbolicEnd}, 
$(\Sigma_{r,\theta})$ converges to $\partial\Cal{E}$ in the $C^0$ sense as $r$ tends to $+\infty$. Since this concerns the convergence of convex functions, it automatically also implies convergence of the spaces of supporting hyperplanes. 
\medskip
\noindent Finally, by Corollary \procref{CurvatureOfLevelSetsII} and the Geometric Maximum Principle (Lemma \procref{LemmaGeomMaxPrinc}), the distance of
$\Sigma_{r,\theta}$ from $\partial_0\Cal{E}$ is at least $R$, where:
$$
\opTanh(R) = \frac{\opTan(\theta-(n-1)\pi/2)}{r}.
$$
\noindent Let $\hat{R}_\theta$ be the maximal value of $R$ which is obtained when $r=\opTan(\theta/n)$:
$$
\opTanh(\hat{R}_\theta) = \frac{\opTan(\theta-(n-1)\pi/2)}{\opTan(\theta/n)}.
$$
\noindent This yields a lower bound for the furthest extent of the foliation for each $\theta$. Since $(\theta-(n-1)\pi/2)(\theta/n)$ converges to $1$ as $\theta$
converges to $n\pi/2$, $\hat{R}_\theta$ converges to $\infty$ as $\theta$ converges to $n\pi/2$ and the result follows.\qed
\medskip
{\bf\noindent Proof of Theorem \procref{TheoremContinuousDependance}:\ }This follows from uniqueness and Lemma \procref{EquivariantDeformation}.\qed
\goodbreak
\newhead{Quasi-Fuchsian Manifolds}
\noindent Quasi-Fuchsian manifolds provide an interesting special case. For all $m$, let $\Bbb{H}^m$ be
$m$-dimensional hyperbolic space. Let $M$ be a compact $n$-dimensional, hyperbolic manifold. We view $\pi_1(M)$ as a subgroup
$\Gamma$ of $\opIsom(\Bbb{H}^n)$.
\headlabel{QuasiFuchsianManifolds}
\medskip
\noindent We denote by $\opRep(\Bbb{H}^n,\Gamma)$ the space of pairs $(\varphi,\alpha)$, where:
\medskip
\myitem{(i)} $\alpha:\Gamma\rightarrow\opIsom(\Bbb{H}^{n+1})$ is a properly discontinous representation of $\Gamma$ in $\opIsom(\Bbb{H}^{n+1})$, and
\medskip
\myitem{(ii)} $\varphi:\partial_\infty\Bbb{H}^n\rightarrow\partial_\infty\Bbb{H}^{n+1}$ is an injective, continuous mapping which is equivariant with respect to $\alpha$.
\medskip
\noindent The set $\opRep(\Bbb{H}^n,\Gamma)$ is a subset of the set of continuous mappings from $\partial_\infty\Bbb{H}^n\munion\Gamma$ into 
$\partial_\infty\Bbb{H}^{n+1}\munion\opIsom(\Bbb{H}^{n+1})$. We furnish this set with the topology of local uniform convergence.
\medskip
\noindent For all $n$, $\Bbb{H}^n$ embeds totally geodesically into $\Bbb{H}^{n+1}$. This induces a homeomorphism $\alpha_0:\opPSO(n,1)\rightarrow\opPSO(n+1,1)$ and an injective
continuous mapping $\varphi_0:\partial_\infty\Bbb{H}^n\rightarrow\partial_\infty\Bbb{H}^{n+1}$ which is equivariant with respect to $\alpha_0$. The connected component of 
$\opRep(\Bbb{H}^n,\Gamma)$ which contains $(\varphi_0,\alpha_0)$ is called the {\bf quasi-Fuchsian} component. The pair $(\varphi, \alpha)$ is then said to be 
{\bf quasi-Fuchsian} if and only if it belongs to the quasi-Fuchsian component.
\medskip
\noindent Let $(\varphi,\alpha)$ be quasi-Fuchsian. Since $\alpha(\Gamma)$ is properly discontinuous, it defines a quotient manifold $\hat{M}_\alpha=\Bbb{H}^{n+1}/\alpha(\Gamma)$. 
When $\alpha=\alpha_0$, we call this manifold the {\bf extension} of $M$.
In the sequel, we identify a quasi-Fuchsian pair and its quotient manifold, and we say that a manifold is
{\bf quasi-Fuchsian} if and only if it is the quotient manifold of a quasi-Fuchsian pair. In this case it may be isotoped
to the extension of a compact, hyperbolic manifold.
\medskip
\noindent Let $(\varphi,\alpha)$ be quasi-Fuchsian. The image of $\partial_\infty\Bbb{H}^n$ under the action of 
$\varphi$ divides $\partial_\infty\Bbb{H}^{n+1}$ into 
two open, simply connected, connected components. The group $\alpha(\Gamma)$ acts properly discontinuously on each of these connected components. The quotient of each 
component is a M\"obius manifold homeomorphic to $M$, and the union of these two quotients forms the ideal boundary of $\hat{M}_\alpha$.
\medskip
\noindent Let $K$ be the convex hull in $\Bbb{H}^{n+1}$ of $\varphi(\partial_\infty\Bbb{H}^n)$. This is the intersection of all closed sets with totally geodesic 
boundary whose ideal boundary does not intersect $\varphi(\partial_\infty\Bbb{H}^n)$. This set is equivariant under the action of $\alpha$ and thus
quotients down to a compact, convex subset of $\hat{M}_\alpha$ which we refer to as the {\bf Nielsen kernel} of $\hat{M}_\alpha$ and which we also denote by $K$. 
Trivally $M\setminus K$ consists of two hyperbolic ends arising from FCSs.
\medskip
\noindent Let $M$ be a quasi-Fuchsian manifold, let $K$ be its Nielsen kernel and let $D$ be the diameter of $K$. Let $\Cal{E}$ be one of the connected components of 
$M\setminus K$. Let $\theta\in[(n-1)\pi/2,n\pi/2[$ be an angle. By Theorem \procref{TheoremExistenceAndUniqueness}, there exists a family $(\Sigma_r)_{r\in]\opTan(\theta/n),\infty[}$ of compact, convex, immersed
hypersurfaces in $\Omega$ such that, for all $r$:
\medskip
\myitem{(i)} $[\Sigma_r]$ is the fundamental class of $\Omega$ and
\medskip
\myitem{(ii)} $R_\theta(\Sigma_r) = r$.
\medskip
\noindent Moreover, this family foliates a neighbourhood of $\partial K\minter\Cal{E}$. We show that this foliation covers the whole of $\Cal{E}$:
\proclaim{Lemma \nextprocno}
\noindent $(\Sigma_r)_{r\in]\opTan(\theta/n),+\infty[}$ foliates the whole of $\Cal{E}$ and $\Sigma_r\rightarrow\partial_\infty\Cal{E}$ in the Hausdorff sense as
$r\rightarrow\opTan(\theta/n)$.
\endproclaim
\proclabel{LemmaUpperBoundInFuchsianMfd}
\proof Let $K_0'$ be the component of $\partial K$ which does not intersect $\Cal{E}$ (i.e. $K_0'$ is the boundary component of $K$ lying on the other side of
$K$ from $\Omega$). For all $d>0$, let $K'_d$ be the level hypersurface in $\Omega\munion K$ at a distance of $d$ from $K_0'$. As in Corollary \procref{CurvatureOfLevelSetsI},
for all $d>0$, the $\theta$-special Lagrangian curvature of $K_d$ is at most $\opTan(\theta/n)/\opTanh(d)$ in the weak sense.
\medskip
\noindent For all $r$, since $\Sigma_r=(S,i_r)$ is compact, there exists a point $p\in S$ such that $d(i_r(p),K_0')$ is minimised. Let $d$ be the distance of $i_r(p)$ from $K_0'$. 
$\Sigma$ is an exterior tangent to $K_d$ at $p$. By the geometric maximum principal:
$$
d(i_r(p),K_0') \geqslant \opArcTanh(r^{-1}\opTan(\theta/n)) - D.
$$
\noindent The result now follows.\qed
\medskip
\noindent The proof of Theorem \procref{TheoremFoliationsInQuasiFuchsian} follows immediately:
\medskip
{\bf\noindent Proof of Theorem \procref{TheoremFoliationsInQuasiFuchsian}:\ } This is the union of Theorem \procref{TheoremExistenceAndUniqueness} and Lemma \procref{LemmaUpperBoundInFuchsianMfd}.\qed
\inappendicestrue
\global\headno=0
\goodbreak
\newhead{Appendix - On a Result of Kamishima}
\noindent An earlier revision of this paper relied on a result of Kamishima (Theorem $B$ of \cite{KamiA}) concerning FCSs whose developing maps are not surjective. We 
discovered that the Kulkarni-Pinkall metric may be used to provide a relatively short proof of this result, which we thus include here.
\medskip
\noindent Let $\Gamma$ be a subgroup of $\opIsom(\Bbb{H}^n)$. The limit set of $\Gamma$, $L(\Gamma)$, is the set of all limit points of sequences of the form 
$(\gamma_n(p))_\ninn$ where $p\in\partial_\infty\Bbb{H}^n$ and $(\gamma_n)_\ninn\in\Gamma$. By definition, this is a closed set. We recall the following important lemma (see, for example \cite{KamiA}):
\proclaim{Lemma \nextprocno, {\bf Chen \& Greenberg, \cite{ChenGreen}}}
\noindent Let $C$ be a closed subset of $\partial_\infty\Bbb{H}^n$ which contains more than one point and is invariant under $\Gamma$, then $L(\Gamma)\subseteq C$.
\endproclaim
\proclabel{ChenGreenberg}
\noindent This yields the following result of Kamishima:
\proclaim{Theorem \nextprocno, {\bf Kamishima, \cite{KamiA}}}
\noindent Let $M$ be a closed conformally flat manifold of dimension at least $3$. If the developing map is not surjective, then it is a covering map.
\endproclaim
\proclabel{ThmKamiCoveringMap}
\proof Let $\tilde{M}$ be the universal cover of $M$, let $\varphi:\tilde{M}\rightarrow\partial_\infty\Bbb{H}^{n+1}$ be 
its developing map and let $\theta:\pi_1(M)\rightarrow\opIsom(\Bbb{H}^{n+1})$ be its holonomy. We consider the two cases where the complement of $\varphi(\tilde{M})$ contains
only one point and where it contains more than one point seperately. Suppose first that $\varphi(\tilde{M})^c$ contains only one point. This point is invariant under the
action of $\Gamma:=\theta(\pi_1(M))$. $\Gamma$ is thus conjugate to a subgroup of the symmetry group of Euclidean space. The result then follows by \cite{Fried}.
Suppose now that $\varphi(\tilde{M})^c$ contains more than one point. Since it is closed and invariant under the action of $\Gamma$, it follows from Lemma \procref{ChenGreenberg} that $L(\Gamma)\subseteq\varphi(\tilde{M})$. In other words, $\varphi(\tilde{M})\subseteq L(\Gamma)^c$. Let $g_{KP}$ be the Kulkarni/Pinkall
metric of $L(\Gamma)^c$ (see \cite{Kulkarni}). since $L(\Gamma)$ contains at least two points, this metric is non-trivial. Moreover, it is complete and invariant under the action of $\Gamma$. Thus $\varphi^*g_{KP}$ is invariant under $\pi_1(M)$. Since $M$ is compact, $\varphi^*g_{KP}$ defines a complete metric over
$\tilde{M}$. $\varphi$ is thus a local isometry between complete manifolds, and the result now follows.\qed
\proclaim{Corollary \nextprocno}
\noindent Let $M$ be a closed conformally flat manifold of dimension at least $3$. If the developing map $\varphi$ is not surjective, then $L(\Gamma)=\partial\varphi(\tilde{M})$.
\endproclaim
\proclabel{LocationOfBoundary}
\goodbreak
\newhead{Bibliography}
{\leftskip = 5ex \parindent = -5ex
\leavevmode\hbox to 4ex{\hfil \cite{AndBarbBegZegh}}\hskip 1ex{Andersson L., Barbot T., B\'e guin F., Zeghib A., Cosmological time versus CMC time in spacetimes of constant curvature}
\medskip
\leavevmode\hbox to 4ex{\hfil \cite{Aubin}}\hskip 1ex{Aubin T., {\sl Nonlinear analysis on manifolds. Monge-Amp\`ere equations}, Die Grund\-lehren der mathematischen Wissenschaften, {\bf 252}, Springer-Verlag, New York, (1982)}
\medskip
\leavevmode\hbox to 4ex{\hfil \cite{CaffNirSpr}}\hskip 1ex{Caffarelli L., Nirenberg L., Spruck J., The Dirichlet problem for nonlinear second-order elliptic equations. I. Monge-Amp\`ere equation. {\sl Comm. Pure Appl. Math} {\bf 37} (1984), no. 3, 369--402} 
\medskip
\leavevmode\hbox to 4ex{\hfil \cite{ChenGreen}}\hskip 1ex{Chen S., Greenberg L., Hyperbolic Spaces, Contribution to {\sl Analysis}, Academic Press, New York, (1974), 49-87}
\medskip
\leavevmode\hbox to 4ex{\hfil \cite{EpsMard}}\hskip 1ex{Epstein D. B. A., Marden, A., Convex hulls in hyperbolic space, a theorem of Sullivan, and measured pleated surfaces, 
In {\sl Fundamentals of hyperbolic geometry: selected expositions}, London Math. Soc. Lecture Note Ser., {\bf 328}, Cambridge Univ. Press, Cambridge, (2006)}
\medskip
\leavevmode\hbox to 4ex{\hfil\cite{Fried}}\hskip 1ex{Fried D., Closed Similarity Manifolds, {\sl Comment. Math. Helvetici} {\bf 55} (1980), 576-582}
\medskip
\leavevmode\hbox to 4ex{\hfil\cite{GuanSpruck}}\hskip 1ex{Guan B., Spruck J., The existence of hypersurfaces of constant Gauss curvature with prescribed boundary, {\sl J. Differential Geom.} {\bf 62} (2002), no. 2, 259--287}
\medskip
\leavevmode\hbox to 4ex{\hfil \cite{HarveyLawson}}\hskip 1ex{Harvey R., Lawson H. B. Jr., Calibrated geometries, {\sl Acta. Math.} {\bf 148} (1982), 47--157}
\medskip
\leavevmode\hbox to 4ex{\hfil \cite{KamiA}}\hskip 1ex{Kamishima T., Conformally Flat Manifolds whose Development Maps are not Surjective, {\sl Trans. Amer. Math. Soc.} {\bf 294} (1986), no. 2, 607-623}
\medskip
\leavevmode\hbox to 4ex{\hfil \cite{KamTan}}\hskip 1ex{Kamishima Y., Tan S., Deformation spaces on geometric structures, In {\sl Aspects of low-dimensional manifolds}, 
Adv. Stud. Pure Math., {\bf 20}, Kinokuniya, Tokyo, (1992)}
\medskip
\leavevmode\hbox to 4ex{\hfil \cite{KapB}}\hskip 1ex{Kapovich M., Deformation spaces of flat conformal structures. Proceedings of the Second Soviet-Japan Joint Symposium of 
Topology (Khabarovsk, 1989), {\sl Questions Answers Gen. Topology} {\bf 8} (1990), no. 1, 253--264}
\medskip
\leavevmode\hbox to 4ex{\hfil \cite{KraSch}}\hskip 1ex{Krasnov K., Schlenker J.M., On the renormalized volume of hyperbolic $3$-manifolds, math.DG/0607081} 
\medskip
\leavevmode\hbox to 4ex{\hfil \cite{Kulkarni}}\hskip 1ex{Kulkarni R.S., Pinkall U., A canonical metric for M\"obius structures and its applications, {\sl Math. Z.} {\bf 216} (1994), no.1, 89--129}
\medskip
\leavevmode\hbox to 4ex{\hfil \cite{LabA}}\hskip 1ex{Labourie F., Un lemme de Morse pour les surfaces convexes, {\sl Invent. Math.} {\bf 141} (2000), 239--297}
\medskip
\leavevmode\hbox to 4ex{\hfil \cite{LabB}}\hskip 1ex{Labourie F., Probl\`eme de Minkowski et surfaces \`a courbure constante dans les vari\'et\'es hyperboliques, {\sl Bull. Soc. Math. Fr.} {\bf 119} (1991), 307-325}
\medskip
\leavevmode\hbox to 4ex{\hfil \cite{MazzPac}}\hskip 1ex{Mazzeo R., Pacard F., Constant curvature foliations in asymptotically hyperbolic spaces, arXiv:0710.2298v2}
\medskip
\leavevmode\hbox to 4ex{\hfil \cite{SmiA}}\hskip 1ex{Smith G., Special Lagrangian curvature, arXiv:math/0506230}
\medskip
\leavevmode\hbox to 4ex{\hfil \cite{SmiNLD}}\hskip 1ex{Smith G., The non-linear Dirichlet Problem in Hadamard Manifolds, arXiv:0908.3590}
\medskip
\leavevmode\hbox to 4ex{\hfil \cite{SmiH}}\hskip 1ex{Smith G., A brief note on foliations of constant Gaussian curvature, arXiv:0802.2202}
\medskip
\leavevmode\hbox to 4ex{\hfil \cite{Thurston}}\hskip 1ex{Thurston W., {\sl Three-dimensional geometry and topology}, Princeton Mathematical Series, {\bf 35}, Princeton University Press,
Princeton, NJ, (1997)}
\medskip
\par
}%
\enddocument